\documentclass{amsproc}
\usepackage{aspmproc}
\usepackage{amssymb}
\usepackage{amsbsy}
\usepackage{amscd}
\usepackage[mathscr]{eucal}
\usepackage{verbatim}
\usepackage{version}
\usepackage{pstricks}
\usepackage{pst-all}
\usepackage{color}

\usepackage[all]{xy}
\usepackage{mathdots}
\makeatletter
\def\cal{\mathcal}
\def\Bbb{\mathbb}
\def\frak{\mathfrak}

\newenvironment{pf*}[1]{\proof[#1]}{\endproof}
\makeatother

\newenvironment{aenume}{%
  \begin{enumerate}%
  }{\end{enumerate}}
\makeatletter
\@ifclasslater{amsproc}{1999/11/24}{}{
Taken from AMSLaTeX ver. 2.0
\renewcommand*\subjclass[2][1991]{%
  \def\@subjclass{#2}%
  \@ifundefined{subjclassname@#1}{%
    \ClassWarning{\@classname}{Unknown edition (#1) of Mathematics
      Subject Classification; using '1991'.}%
  }{%
    \@xp\let\@xp\subjclassname\csname subjclassname@#1\endcsname
  }%
}
\renewcommand{\subjclassname}{%
  \textup{1991} Mathematics Subject Classification}
\@xp\let\csname subjclassname@1991\endcsname \subjclassname
\@namedef{subjclassname@2000}{%
  \textup{2000} Mathematics Subject Classification}
}
\makeatother
\newenvironment{NB}{
\color{red}{\bf NB}. \footnotesize 
}{}
\excludeversion{NB}
\newenvironment{NB2}{
\color{blue}{\bf NB}. \footnotesize
}{}

\newtheorem{Theorem}[equation]{Theorem}
\newtheorem{Corollary}[equation]{Corollary}
\newtheorem{Lemma}[equation]{Lemma}
\newtheorem{Proposition}[equation]{Proposition}

\theoremstyle{definition}
\newtheorem{Definition}[equation]{Definition}

\theoremstyle{remark}
\newtheorem{Remark}[equation]{Remark}
\newtheorem*{Claim}{Claim}

\numberwithin{equation}{section}

\newcommand{\thmref}[1]{Theorem~\ref{#1}}
\newcommand{\secref}[1]{\S\ref{#1}}
\newcommand{\lemref}[1]{Lem\-ma~\ref{#1}}
\newcommand{\propref}[1]{Proposition~\ref{#1}}
\newcommand{\corref}[1]{Corollary~\ref{#1}}
\newcommand{\subsecref}[1]{\S\ref{#1}}

\newcommand{\defref}[1]{Definition~\ref{#1}}
\newcommand{\remref}[1]{Remark~\ref{#1}}
\newcommand{\lsp}[2]{\,{}^{#1}{#2}}
\newcommand{\defeq}{{:=}}
\newcommand{\C}{{\Bbb C}}
\newcommand{\Z}{{\Bbb Z}}
\newcommand{\Q}{{\Bbb Q}}
\newcommand{\R}{{\Bbb R}}
\newcommand{\proj}{{\Bbb P}}

\newcommand{\SL}{\operatorname{\rm SL}}

\newcommand{\GL}{\operatorname{GL}}

\newcommand{\algsl}{\operatorname{\frak{sl}}} 

\newcommand{\End}{\operatorname{End}}
\newcommand{\Hom}{\operatorname{Hom}}
\newcommand{\Ext}{\operatorname{Ext}}
\newcommand{\Ker}{\operatorname{Ker}}
\newcommand{\Coker}{\operatorname{Coker}}
\newcommand{\Ima}{\operatorname{Im}}

\newcommand{\rank}{\operatorname{rank}}

\newcommand{\tr}{\operatorname{tr}}

\newcommand{\id}{\operatorname{id}}

\newcommand{\linf}{{\ell_\infty}}
\newcommand{\shfO}{\mathcal O}
\newcommand{\dslash}{/\!\!/} 
\newcommand{\bp}{{{\widehat\proj}^2}}
\newcommand{\bM}{{\widehat M}}

\newcommand{\ch}{\operatorname{ch}}
\newcommand{\Wedge}{{\textstyle \bigwedge}}

\newcommand{\hT}{\widetilde T}

\def\<{\langle}
\def\>{\rangle}

\newcommand{\bMreg}{\widehat M^{\operatorname{reg}}_0}

\newcommand{\codim}{\mathop{\text{\rm codim}}\nolimits}

\newcommand{\tzeta}{{}^0\zeta}
\newcommand{\Qcal}{\mathcal Q}

\makeatletter
\newcommand{\vechatom}{
    {\Vec{\omega}}
    \,\smash[b]{\hbox{\lower2\ex@\hbox{$\m@th\hat{\null}$}}}
}
\makeatother

\usepackage[dvips,bookmarks=false]{hyperref}

\begin{document}
\title[Perverse coherent sheaves on blow-up. I]
{Perverse coherent sheaves on blow-up. I.
\\ a quiver description
}

\author[H.~Nakajima and K.~Yoshioka]{Hiraku Nakajima and K\={o}ta Yoshioka}
\address{
{\rm Hiraku Nakajima}\\
Department of Mathematics\\
Kyoto University\\
Kyoto 606-8502\\
Japan
}
\email{nakajima@math.kyoto-u.ac.jp}
\address{
{\rm K\={o}ta Yoshioka}\\
Department of Mathematics, Faculty of Science\\
Kobe University\\
Kobe 657-8501\\
Japan
}
\email{yoshioka@math.kobe-u.ac.jp}


\maketitle
\section*{Introduction}

This is the first of two papers studying moduli spaces of a certain class of
coherent sheaves, which we call {\it stable perverse coherent
  sheaves}, on the blowup of a projective surface. They are used to
relate usual moduli spaces of stable sheaves on a surface and its
blowup.

Let us give the definition for general case though we will consider
only framed sheaves on the blowup $\bp$ of the projective plane
$\proj^2$ in this paper.
Let $p\colon \widehat X\to X$ be the blowup of a smooth projective
surface $X$ at a point. Let $C$ be the exceptional divisor. 
A {\it stable perverse coherent sheaf\/} $E$ on $\widehat X$,
with respect to an ample line bundle $H$ on $X$,
is
\begin{enumerate}
\item
$E$ is a coherent sheaf on $\widehat X$,
\item
$\Hom(E(-C),{\cal O}_C)=0$,
\item
$\Hom(\shfO_C,E) = 0$,
\item $p_*E$ is slope stable with respect to $H$.%
\begin{NB}
\item
$R^1 p_*(E)=0$,

In the first version of Kota's note, this condition was imposed. But
it turns out that this condition automatically follows from the other
conditions. This is also observed in \subsecref{subsec:higherdirect}.
\end{NB}
\end{enumerate}
We will consider the moduli spaces of coherent sheaves $E$ on $\widehat X$
such that $E(-mC)$ is stable perverse coherent for $m\ge 0$. This
depends on $m$.

Let us first explain how we find the definition of stable perverse
coherent sheaves.
Our definition comes from two sources. The first one is Bridgeland's
paper \cite{Br:4}, from which we take the name ``perverse coherent''
sheaves. He introduced a new $t$-structure on the derived category
${\mathrm D}^b(\operatorname{Coh} Y)$ of coherent sheaves on $Y$ for a
birational morphism $f\colon Y\to X$ such that (1) ${\mathbf
  R}f_*(\shfO_Y) = \shfO_X$ and (2) $\dim f^{-1}(x) \le 1$ for any
$x\in X$. An object $E\in {\mathrm D}^b(\operatorname{Coh} Y)$ is
called {\it perverse coherent\/} if it is in the heart
$\operatorname{Per}(Y/X)$ of the $t$-structure.
Then he considered the moduli spaces of perverse ideal sheaves which
are subobjects of the structure sheaf $\shfO_Y$ in the heart. His
interest was in the case when $\dim Y = 3$ (and to describe the flop
of $f\colon Y\to X$ as moduli space of perverse ideal sheaves). But
his definition of the $t$-structure makes sense for the blowup
$p\colon \widehat X\to X$ of a smooth quasi-projective surface $X$ at a
point.
Then we get conditions similar as above, but we further impose the
stability condition (4), then it implies that $E$ is a coherent sheaf,
not an object in the derived category. (This will be explained in the
subsequent paper.)

The second source, which is explained in detail in this paper, is the
quiver description of the framed moduli spaces of locally free sheaves
on the blowup $\bp$ of the projective plane due to King~\cite{King},
which is based on an earlier work by Buchdahl~\cite{Bu}.
Here the framing is a trivialization of the sheaf at the line at infinity.
He described the relevant moduli spaces as GIT (geometric invariant
theory) quotients of representations of a certain quiver with a
relation. (See \secref{sec:ADHM} for a precise statement.) The groups
taking quotients are the products of two copies of general linear
groups of possibly different sizes. A particular linearization was
used to define the GIT quotients in the original paper, but we can
consider more general linearization. Then we will show that more
general GIT quotients parametrize framed perverse coherent sheaves
(after twisting by $\shfO(-mC)$),
where we impose conditions (1)$\sim$(3) above, but we replace (4) by 
\begin{enumerate}
\item[(4)'] $E$ is torsion free outside $C$.
\end{enumerate}
This is because the stability is implicit in the existence of the
framing. Note also that (3) is equivalent to the torsion freeness of
$p_*E$ at the point which we blowup (see the proof of
\lemref{lem:cohvanish}). Therefore (3) and (4)' are combined to the
condition that $p_*E$ is torsion free. But we keep them separated, as
(3) is be altered by a twist by $\shfO(C)$ while (4)' is not.

\begin{figure}[htbp]
\def\JPicScale{.8}
  \centering
\ifx\JPicScale\undefined\def\JPicScale{1}\fi
\psset{unit=\JPicScale mm}
\psset{linewidth=0.3,dotsep=1,hatchwidth=0.3,hatchsep=1.5,shadowsize=1,dimen=middle}
\psset{dotsize=0.7 2.5,dotscale=1 1,fillcolor=black}
\psset{arrowsize=1 2,arrowlength=1,arrowinset=0.25,tbarsize=0.7 5,bracketlength=0.15,rbracketlength=0.15}
\begin{pspicture}(0,0)(90,65)
\psline(30,10)(90,10)
\psline[linewidth=0.4,linestyle=dotted](90,0)(30,60)
\rput(80,65){$\zeta_0 = 0$}
\rput[r](25,10){$\scriptscriptstyle\zeta_1 = 0$}
\rput[l](90,5){$\zeta_0+\zeta_1 = 0$}
\newrgbcolor{userFillColour}{0.8 0.8 0.8}
\pspolygon[linestyle=none,fillcolor=userFillColour,fillstyle=solid](30,60)
(80,10)
(30,53.75)(30,60)
\pspolygon[linestyle=none,fillcolor=lightgray,fillstyle=solid](30,10)(80,10)(80,0)(30,0)
\psline[linewidth=0.4,linestyle=dotted](80,0)(80,60)
\psline(30,35)(80,10)
\psline(30,43.33)(80,10)
\psline(30,47.5)(80,10)
\psline(30,50)(80,10)
\psline(30,51.67)(80,10)
\psbezier(30,52.87)(30,52.86)(80,10)(80,10)
\psline(30,53.75)(80,10)
\rput[r](25,35){$\scriptscriptstyle\zeta_0+2\zeta_1 = 0$}
\rput[r](25,43.33){$\scriptscriptstyle 2\zeta_0+3\zeta_1=0$}
\rput[r](25,47.5){$\scriptscriptstyle 3\zeta_0+4\zeta_1=0$}
\rput[r](25,50){$\scriptscriptstyle 4\zeta_0+5\zeta_1=0$}
\rput[l](85,20){${}^0\zeta$}
\rput{0}(66.5,5.5){\psellipse[fillstyle=solid](0,0)(0.5,-0.5)}
\psline[linewidth=0.15]{<-}(67,6)(85,20)
\rput{0}(35.5,52.5){\psellipse[fillstyle=solid](0,0)(0.5,-0.5)}
\psline[linewidth=0.15]{<-}(36,53)(45,60)
\rput[l](45,60){$^\infty\zeta$}
\rput[r](15,55){$\scriptscriptstyle\vdots$}
\end{pspicture}
\caption{Chamber structure}
\label{fig:zeta1}
\end{figure}
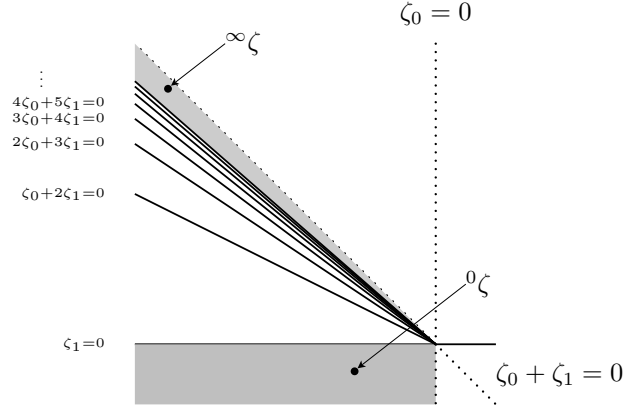

From either view point, we have a chamber structure on the stability
parameter. In the first view point, it is the choice of $m$. In the
second, we have a parameter $\zeta = (\zeta_0,\zeta_1)\in \R^2$. We
are interested in the domain $\zeta_0+\zeta_1 < 0$, $\zeta_0 < 0$ in
this paper. This domain is separated by walls defined by $m \zeta_0 +
(m+1)\zeta_1 = 0$ for $0\le m \le N$ for some $N$ depending on Chern
classes of sheaves. (See Figure~\ref{fig:zeta1}.)

The wall $m \zeta_0 + (m+1)\zeta_1 = 0$ corresponds to the sheaf
$\shfO_C(-m-1)$. When we cross the wall, sheaves given by extensions
$0\to E'\to E^+ \to \shfO_C(-m-1)^{\oplus p}\to 0$ are replaced by
extensions in the opposite way $0\to \shfO_C(-m-1)^{\oplus p}\to E^-
\to E'\to 0$. This is very similar to the wall-crossing phenomenon
appearing when we change polarization in the moduli spaces of sheaves.
We study this in detail in the subsequent paper.

There are two distinguished chambers in our domain, the gray regions
in the figure. The chamber adjacent to the plane $\zeta_0+\zeta_1 = 0$
(containing the parameter ${}^\infty\zeta$ in the figure) gives 
moduli spaces of usual framed sheaves (and $(p^*H-\varepsilon
C)$-stable sheaves for an ample line bundle $H$ on $X$ and
sufficiently small $\varepsilon$ for the projective case).
King's stability condition, corresponding to framed {\it locally
  free\/} sheaves, lies on the plane $\zeta_0+\zeta_1 = 0$.
On the other hand, the chamber $\zeta_0 < 0$, $\zeta_1 < 0$
(containing the parameter ${}^0\zeta$ in the figure) gives moduli
spaces of framed perverse coherent sheaves. We will show that these
moduli spaces are isomorphic to those of framed sheaves on $\proj^2$,
instead of $\bp$ (and $H$-stable sheaves on $X$ for the projective
case) if the first Chern class vanishes in \subsecref{subsec:another2}.
(We have a similar statement for general first Chern classes.)
The simplest (nontrivial) case of this statement says that the framed
moduli space of perverse ideal sheaves on $\bp$ of colength $1$ is
isomorphic to $\C^2$, which is the framed moduli space of ideal
sheaves on $\proj^2$.
This statement is natural in view of Bridgeland's result~\cite{Br:4}.
In the $3$-dimensional case $f\colon Y\to X$, the moduli space of
perverse ideal sheaves of colength $1$ is the flop of $Y$.
However the derived categories
${\mathrm D}^b(\operatorname{Coh} X)$ and
${\mathrm D}^b(\operatorname{Coh} \widehat X)$,
are {\it not\/} equivalent like the flop case, so our analogy breaks
at this point.

Our result should have generalization to various other situations. For
example, we can consider a similar problem for the minimal resolution
$p\colon Y \to X$ where $X$ has a quotient singularity of the form
$\C^2/\Gamma$ for a finite subgroup $\Gamma$ of $\GL_2(\C)$.
In particular, when $\Gamma$ is a subgroup of $\SL_2(\C)$, there is a
quiver description of framed sheaves on the minimal resolution of
$\C^2/\Gamma$ (\cite{KN,Na:ADHM}). We have a chamber structure on the
stability condition, where the space of parameters is identified with
the Cartan subalgebra of the affine Lie algebra corresponding to
$\Gamma$, and walls are root hyperplanes. Moduli spaces parametrize
perverse coherent sheaves twisted by line bundles, as in this paper.
We again have two distinguished chambers. As in this paper,
${}^\infty\zeta$ gives the moduli of framed torsion free sheaves, while
${}^0\zeta$ gives the moduli of framed $\Gamma$-equivariant torsion
free sheaves on $\C^2$.
These moduli spaces are studied, in the name of {\it quiver
  varieties\/}, by the first named author \cite{Na:1994,Na:1998}.
When sheaves have rank $1$, moduli spaces of perverse coherent sheaves
were constructed and studied by Toda \cite{Toda}.
However there is a {\it crucial difference\/} in the wall-crossing. In
this case, all moduli spaces are diffeomorphic to each other, as
proved in \cite{Na:1994}. This is not true in our case, as framed
moduli spaces of sheaves on $\bp$ and $\proj^2$ have different Betti
numbers. (See \cite{NY2}.) 
If we see Figure~\ref{fig:zeta1}, it looks like the Cartan
subalgebra for the affine $\algsl_2$ with root hyperplanes. However if
we draw all walls without imposing the condition $\zeta_0+\zeta_1 <
0$, $\zeta_0 < 0$, then it becomes different. (See
\remref{rem:allW=0}.)
It should be also noticed that $\shfO_C(-m-1)$ is an exceptional
object (i.e., has no higher self-extensions) while the corresponding
objects in the $\SL_2(\C)$-setting are spherical objects.
If $\Gamma$ is not contained in $\SL_2(\C)$, there should be
difference from $\SL_2(\C)$ cases, but we do not understand the
picture so far.

It is also interesting to study the problem in Bridgeland's
$3$-dimen\-sional case \cite{Br:4}, and also in more general situation
\cite{Br:2007}, when we could expect an explicit description of all
{\it walls\/} as in this paper.

Our motivation to this research comes from the authors' study on
the `instanton counting'. (See \cite{part1,NY2,NY3} and the
references therein.)
An understanding of relations between moduli spaces of sheaves on $X$
and $\widehat X$ was one of the most essential ingredient there. An
application of our result to the instanton counting will be discussed
in a separate paper.

\subsection*{Acknowledgments}

The first named author is supported by the Grant-in-aid for Scientific
Research (No.\ 19340006), JSPS. A part of this paper was written while
the first named author was visiting the Institute for Advanced Study
with supports by the Ministry of Education, Japan and the Friends of
the Institute.

\subsection*{Notations}

Let $[z_0:z_1:z_2]$ be the homogeneous coordinates on $\proj^2$ and
$\linf = \{ z_0 = 0\}$ the line at infinity.
Let $p\colon \bp\to \proj^2$ be the blowup of $\proj^2$ at $[1:0:0]$.
Then $\bp$ is the closed
subvariety of $\proj^2\times\proj^1$ defined by
\begin{equation*}
  \left\{ ([z_0:z_1:z_2],[z:w])\in\proj^2\times\proj^1 \mid
    z_1 w = z_2 z \right\},
\end{equation*}
where the map $p$ is the projection to the first factor.  We denote
$p^{-1}(\linf)$ also by $\linf$ for brevity. Let $C$ denote the
exceptional divisor given by $z_1 = z_2 = 0$. Let $\shfO$ denote the
structure sheaf of $\bp$, $\shfO(C)$ the line bundle associated with
the divisor $C$, and $\shfO(mC)$ its $m^{\mathrm{th}}$ tensor product
$\shfO(C)^{\otimes m}$ when $m > 0$, $\left(\shfO(C)^{\otimes
    -m}\right)^\vee$ if $m < 0$, and $\shfO$ if $m=0$. And we use the
similar notion $\shfO(mC+n\linf)$ for tensor products of $\shfO(mC)$
and tensor powers of the line bundle corresponding to $\linf$ or its
dual.

The structure sheaf of the exceptional divisor $C$ is denoted by
$\shfO_C$. If we twist it by the line bundle $\shfO_{\proj^1}(n)$ over
$C \cong\proj^1$, we denote the resulted sheaf by $\shfO_C(n)$. Since
$C$ has the self-intersection number $-1$, we have
$\shfO_C\otimes \shfO(C) = \shfO_C(-1)$.

\begin{NB}
I do not decide where to put the following subsection.

I decided to put this into the proof of \lemref{lem:cohvanish}.

\subsection{}

Following Kota's note, we say a framed sheaf $(E,\Phi)$ is {\it stable
  perverse coherent\/} if and only if the following holds:
\begin{enumerate}
\item
$E$ is a coherent sheaf on $\bp$,
\item
$\Hom(E,{\cal O}_C(-1))=0$,
\item
$p_*(E)$ is a torsion free sheaf on $\proj^2$,
\item $E$ is locally free of rank $r$ in a neighborhood of $\linf$ and
  $\Phi\colon E|_{\linf}\to \shfO_{\linf}^{\oplus r}$ is an
  isomorphism called a `framing'.
\item
$R^1 p_*(E)=0$,

In the first version of Kota's note, this condition was imposed. But
it turns out that this condition automatically follows from the other
conditions. This is also observed in \subsecref{subsec:higherdirect}.
\end{enumerate}

We see that the condition (3) is equivalent to 
\begin{enumerate}
\item[(3)'] $\Hom(\shfO_C,E) = 0$.
\end{enumerate}
In fact, let $\C_0$ be the skyscraper sheaf at $[1:0:0]$ of $\proj^2$.
Then we have $\Hom(\shfO_C,E) = \Hom(p^*\C_0, E) = \Hom(\C_0,p_*E)$.
In this note, (3)' is more natural from \propref{prop:ss}.

We give one more supplement to Kota's note:

If $R^1 p_*(E) = 0$, then $R^1p_* E(-mC) = 0$ for $m > 0$. (This is
easier version of Lemma~1.8 in his note.) In fact, we have $p^*
p_*(\shfO(-mC)) \twoheadrightarrow \shfO(-mC)$ as $m > 0$. We further
take a vector bundle $V$ such that $V\twoheadrightarrow
p_*(\shfO(-mC))$. Then we have $E\otimes p^* V \twoheadrightarrow
E(-mC)$. We take the higher direct image and use the projection
formula to get
\begin{equation*}
    0 = R^1p_*(E) \otimes V\to R^1p_*(E(-mC))
    \to 0,
\end{equation*}
where the final term is $R^2$ of something, but it vanishes as the
fiber of $p$ has dimension at most $1$.
\end{NB}

\section{ADHM description -- Main result}\label{sec:ADHM}

\subsection{Preliminary}
Let $M(r,n)$ be the framed moduli space of torsion free sheaves $E$ on
$\proj^2 = \C^2\cup\linf$ with $\rank E = r > 0$, $c_2(E)[\proj^2] = n
\ge 0$, where $E$ is assumed to be locally free along $\linf$ and the
framing is a trivialization on the line at infinity $\Phi\colon
E|_{\linf}\to \shfO_{\linf}^{\oplus r}$.
The ADHM description identifies $M(r,n)$ with the space of the
following data $X = (B_1, B_2, i, j)$ defined for vector spaces $V$,
$W$ with $\dim V = n$, $\dim W = r$ modulo the action of $\GL(V)$:
\begin{itemize}
\item $B_1, B_2\in \End(V)$, $i\in\Hom(W,V)$, $j\in\Hom(V,W)$,
\item $[B_1, B_2] + ij = 0$,
\item (stability condition) \mbox{} \\
  a subspace $T\subset V$ with $B_\alpha(T)\subset T$
  ($\alpha=1,2$), $\Ima i\subset T$ must be $T = V$.
\end{itemize}

The quotient space of data by the action by $\GL(V)$ can be understand
as the GIT (geometric invariant theory) quotient with respect to the
trivial line bundle with the $\GL(V)$-action given by
$\det\colon\GL(V)\to\C^*$. The above stability condition is nothing
but one in GIT.
If we do not impose the last condition, but take the affine
algebro-geometric quotient $\dslash\GL(V)$, we get the Uhlenbeck
(partial) compactification 
\begin{equation*}
  M_0(r,n) = \bigsqcup_{m=0}^n M^{\operatorname{reg}}_0(r,n-m)\times S^m\C^2,
\end{equation*}
where $M^{\operatorname{reg}}_0(r,n-m)$ is the framed moduli space of
locally free sheaves, which is identified with the data
$X=(B_1,B_2,i,j)$ satisfying the above conditions together with
\begin{itemize}
\item (co-stability condition) \mbox{} \\
  a subspace $S\subset V$ with $B_\alpha(S)\subset S$
  ($\alpha=1,2$), $\Ker j\supset S$ must be $S = 0$.
\end{itemize}
(See \cite[Chap.\ 2]{Lecture} for detail.)

There is an ADHM type description of a framed moduli space
\allowbreak
$\bMreg(r,k,n)$ of locally free sheaves $E$ on the blowup $\bp$ for
$\rank E = r > 0$, $(c_1(E),[C]) = -k$, $(c_2(E) -
(r-1)c_1(E)^2/(2r),[\bp]) = n \ge 0$ due to King \cite{King}. We take
vector spaces $V_0$, $V_1$, $W$ with
\begin{equation*}
  r = \dim W = r, \quad k = -\dim V_0 + \dim V_1, \quad
 n + \frac{k^2}{2r} = \frac12(\dim V_0+\dim V_1).
\end{equation*}

\begin{Theorem}[\protect{\cite{King}}]\label{thm:King}
The framed moduli space $\bMreg(r,k,n)$ is bijective to
the space of the following data $X = (B_1,B_2,d,i,j)$
modulo the action of $\GL(V_0)\times \GL(V_1)$:
\begin{itemize}
\item $B_1, B_2\in \Hom(V_1, V_0)$, $d\in \Hom(V_0,V_1)$, $i\in
  \Hom(W,V_0)$, $j\in \Hom(V_1,W)$,
\begin{equation*}
\xymatrix{
V_0 \ar@<-1ex>[rr]_{d} && \ar@<-1ex>[ll]_{B_1,B_2} \ar[ll] \ar[ld]^j V_1 \\
& \ar[lu]^i W & \\
}
\end{equation*}
\item $\mu(B_1,B_2,d,i,j) = B_1 d B_2 - B_2 d B_1 + ij = 0$,
\item \textup(stability condition\textup)
\begin{description}
\item[(S1)] for subspaces $S_0\subset V_0$, $S_1\subset V_1$ such that
$B_\alpha(S_1)\subset S_0$ ($\alpha=1,2$),
$d(S_0)\subset S_1$, $\Ker j\supset S_1$, we have $\dim S_0 > \dim
S_1$ or $S_0 =  S_1 = 0$.
\item[(S2)] for subspaces $T_0\subset V_0$,
$T_1\subset V_1$ such that $B_\alpha(T_1)\subset T_0$ ($\alpha=1,2$),
$d(T_0)\subset T_1$, $\Ima i\subset T_0$, we have $\codim T_1 > \codim
T_0$ or $(T_0, T_1) = (V_0, V_1)$,
\end{description}
\end{itemize}
\end{Theorem}

\begin{NB}
In \cite{King} it is assumed that $c_1(E) = 0$. But it is unnecessary.
\end{NB}

\begin{NB}
The torus action $\hT\curvearrowright\bMreg(r,k,n)$ is given by
\begin{equation}\label{eq:action}
   [(B_1, B_2, d, i, j)] \longmapsto
   [(t_1 B_1, t_2 B_2, d, ie^{-1}, t_1 t_2 e j)]
\end{equation}
\end{NB}

The quotient space can be again understood as the GIT quotient with
respect to the trivial line bundle with the
$\GL(V_0)\times\GL(V_1)$-action given by
$\det_{\GL(V_1)}/\det_{\GL(V_0)}$.
See \cite[(3.1)]{King-mod}, \cite[(6.2)]{Na-var-quiver} for the
equivalence between the usual GIT-stability with respect to the line
bundle and the above stability condition in the context of quiver
representations.

\begin{NB}
  The corresponding line bundle on the quotient space is nothing but
  the determinant line bundle $L$ such that $c_1(L) = \mu(C)$. (Need
  check for the sign.)
\end{NB}

A framed locally-free sheaf $(E,\Phi)$ is constructed from $X$ as
follows. Let us consider coordinates $z_0, z_1, z_2$ (resp.\ $z$, $w$)
as sections of $\shfO(\linf)$ (resp.\ $\shfO(-C+\linf)$). The ratio $s
\defeq z_1/z = z_2/w$ is a section of $\shfO(C)$ which vanishes on
$C$. Then the locally free sheaf $E$ corresponding to
$(B_1,B_2,d,i,j)$ is the middle cohomology of the complex
\begin{equation}\label{eq:cpx}
\begin{CD}
\begin{matrix}
   V_0 \otimes \shfO(C-\linf)
\\
   \oplus
\\ V_1\otimes\shfO(-\linf)
\end{matrix}
@>{\alpha}>>
\begin{matrix}
   \C^2\otimes V_0\otimes\shfO \\ \oplus \\ \C^2\otimes V_1\otimes\shfO \\
   \oplus \\ W\otimes\shfO
\end{matrix}
@>{\beta}>>
\begin{matrix}
  V_0\otimes \shfO(\linf)
\\
  \oplus
\\
  V_1\otimes\shfO(-C+\linf)
\end{matrix},
\end{CD}
\end{equation}
with
\begin{equation*}
   \alpha = 
   \begin{bmatrix}
     z & z_0 B_1           \\
     w & z_0 B_2           \\
     0 & z_1 - z_0 d B_1   \\
     0 & z_2 - z_0 d B_2   \\
     0 & z_0 j
   \end{bmatrix},
\qquad
   \beta =
   \begin{bmatrix}
    z_2 & - z_1 & B_2 z_0 & -B_1 z_0 & i z_0\\
    dw  & - d z & w       & -z       & 0
   \end{bmatrix}.
\end{equation*}
The equation $\mu(B_1,B_2,d,i,j) = B_1 d B_2 - B_2 d B_1 + ij = 0$ is
equivalent to $\beta\alpha = 0$. The stability condition (S2) is
equivalent to the surjectivity of $\beta$. The stability condition
(S1) is equivalent to that $\alpha$ is injective and the image is a
subbundle. (See \cite[4.1.3]{King} and also \lemref{lem:King}.)
Therefore $E = \Ker\beta/\Ima\alpha$ is a locally free sheaf as
expected.%
\begin{NB}
\begin{equation*}
\begin{gathered}
  - k = (c_1(E),[C]) = - (c_1(\shfO(-\dim V_1 C + \dim V_0 C), [C])
  = \dim V_0 - \dim V_1,
\\
\begin{aligned}[t]
  n & = - \int_{\bp} \ch_2(E) + \frac1{2r} \int_{\bp} c_1(E)^2
  = \left(\dim V_0 + \dim V_1\right) \left(\int_{\bp} \ch_2(\shfO(\linf))
    + \int_{\bp} \ch_2(\shfO(C-\linf))\right)
  - \frac{k^2}{2r}
\\
  & = \frac12 \left(\dim V_0 + \dim V_1\right)  - \frac{k^2}{2r}.
\end{aligned}
\end{gathered}
\end{equation*}
\end{NB}

Let ${\widehat\C}^2$ be the blowup of $\C^2$ at the origin. The
Uhlenbeck (partial) compactification $\bM_0(r,k,n)$, which is defined
as a set by
\begin{equation*}
   \bM_0(r,k,n) =
   \bigsqcup_{m=0}^{[n]}
   \bMreg(r,k,n-m)\times S^m{\widehat\C}^2
\end{equation*}
is identified with the GIT quotient associated with the line bundle
associated with the above homomorphism. Namely we replace $>$ by $\ge$
in the two inequalities in the stability condition, and then divide by
the $S$-equivalence relation (see \cite[2.3]{King-mod},
\cite[2.3]{Na-var-quiver}).

We can generalize this description to the case when $E =
\Ker\beta/\Ima\alpha$ is only assumed to be {\it torsion free}. 
The corresponding generalization for $\proj^2$ was given in
\cite[Chap.~2]{Lecture}, and it is given by allowing $\alpha$ to be
injective possibly except finitely many points. This is equivalent to
that $\alpha$ is injective as a sheaf homomorphism and
$\Ker\beta/\Ima\alpha$ is torsion free.
The same works in the blowup case. (The proof of this statement, as
well as that of \thmref{thm:King}, will be given in this paper in more
general setting.)
The corresponding equivariant line bundle for the GIT stability is
given by the homomorphism
\(
   \det_{\GL(V_1)}^N / \det_{\GL(V_0)}^{N+1}
\)
for sufficiently large $N$. 
Then the condition (S2) is unchanged, but we replace $>$ in (S1) by
$\ge$.
Let us denote this condition by (S1)'.

\begin{NB}
The torus action on $\bM(r,k,n)$ is given by the same formula as
\eqref{eq:action}.
\end{NB}

This observation naturally leads to consider more general stability
conditions. Let $\zeta = (\zeta_0,\zeta_1)\in\R^2$.

\begin{Definition}\label{def:stable}
Suppose $\dim W = r\neq 0$ as above.
We say $X = (B_1,B_2,d,i,j)$ is {\it $\zeta$-semistable\/} if
\begin{enumerate}
\item for subspaces $S_0\subset V_0$, $S_1\subset V_1$ such that
$B_\alpha(S_1)\subset S_0$ ($\alpha=1,2$),
$d(S_0)\subset S_1$, $\Ker j\supset S_1$, we have 
\(
    \zeta_0 \dim S_0 + \zeta_1 \dim S_1 \le 0.
\)
\item for subspaces $T_0\subset V_0$,
$T_1\subset V_1$ such that $B_\alpha(T_1)\subset T_0$ ($\alpha=1,2$),
$d(T_0)\subset T_1$, $\Ima i\subset T_0$, we have 
\( 
   \zeta_0 \codim T_0 + \zeta_1 \codim T_1 \ge 0.
\)
\end{enumerate}
We say $X$ is {\it $\zeta$-stable\/} if the inequalities are strict
unless $(S_0,S_1) = (0,0)$ and $(T_0,T_1) = (V_0,V_1)$ respectively.
\end{Definition}

Thus the condition~(S1),(S2) is equivalent to $\zeta$-stability with
$\zeta_0 + \zeta_1 = 0$, $\zeta_0 < 0$. In Figure~\ref{fig:zeta1} it
is a parameter on the line $\zeta_0+\zeta_1 = 0$ which is on the
boundary on the domain $\zeta_0 + \zeta_1 < 0$, $\zeta_0 <0$ we are
interested in this paper.

This stability condition come from a $\Q$-line bundle with an
$\GL(V_0)\times\GL(V_1)$-action if $(\zeta_0,\zeta_1)\in \Q^2$, i.e.,
the trivial line bundle with the equivariant structure given by the
$\Q$-homomorphism
\(
   \det_{\GL(V_0)}^{\zeta_0} \det_{\GL(V_1)}^{\zeta_1}.
\)
As we can move the parameter in a chamber (see
\subsecref{subsec:walldef}) without changing the stability condition,
the condition $(\zeta_0,\zeta_1)\in\Q^2$ is not essential.

We consider two types of quotient spaces:
\begin{equation*}
\begin{split}
   & \bM_\zeta^{\mathrm{s}}(r,k,n) \defeq
    \left\{ (B_1,B_2,d,i,j) \in \mu^{-1}(0) \left|\,
    \text{$\zeta$-stable}      
  \right\}\right/ \GL(V_0)\times \GL(V_1),
\\
   & \bM_\zeta^{\mathrm{ss}}(r,k,n) \defeq
    \left\{ (B_1,B_2,d,i,j) \in \mu^{-1}(0) \left|\,
    \text{$\zeta$-semistable}      
  \right\}\right/ \!\sim,
  \end{split}
\end{equation*}
where $\sim$ denotes the $S$-equivalence relation.
(Again, see \cite[2.3]{King-mod}, \cite[2.3]{Na-var-quiver} more detail.)

\subsection{The statement}
We fix $m\in\Z_{\ge 0}$ and assume further
\begin{equation}\label{eq:overallass}
   \zeta_0 < 0, \qquad
   0 < - \left(m\zeta_0 + (m+1)\zeta_1\right) \ll 1.
\end{equation}
In \subsecref{subsec:walldef} we will prove that this condition
specifies a chamber, where the $\zeta$-stability and
$\zeta$-semistability is equivalent. In Figure~\ref{fig:zeta1}, the
point is in the chamber just below the line $m\zeta_0 + (m+1)\zeta_1 =
0$.

The following is the main result of this paper.

\begin{Theorem}\label{thm:main}
  The space $\bM_\zeta^{\mathrm{s}}(r,k,n) =
  \bM_\zeta^{\mathrm{ss}}(r,k,n)$ is bijective to the space of
  isomorphism classes of framed sheaves $(E,\Phi)$ on $\bp$ such that
  $E(-mC)$ is perverse coherent, i.e., it satisfies
\begin{enumerate}
\item
$\Hom(E(-mC),{\cal O}_C(-1))=0$,
\item
$\Hom(\shfO_C,E(-mC)) = 0$,
\item
$E(-mC)$ is torsion free outside $C$.
\end{enumerate}
\end{Theorem}

A construction of moduli spaces of stable perverse coherent sheaves
will be given in the subsequent paper. It can be adapted to a
construction of the framed moduli spaces.
However our proof also gives an universal family on the space of
$\zeta$-stable data, so our result can be also considered as a
(different) construction of the fine moduli spaces of framed stable
perverse coherent sheaves.

\subsection{Morphisms between moduli spaces}

Let us consider the affine algebro-geometric quotient
$\mu^{-1}(0)\dslash \GL(V_0)\allowbreak
\times \GL(V_1)$. We have a morphism
$\mu^{-1}(0)\dslash \GL(V_0)\times \GL(V_1)\to M_0(r,n)$
($n = \min(\dim V_0,\dim V_1))$ given by
\begin{equation*}
    [(B_1,B_2,d,i,j)]\longmapsto [(dB_1, dB_2, di, j)] \text{ or }
    [(B_1 d, B_2 d, i, jd)].
\end{equation*}
The coordinate ring of $\mu^{-1}(0)\dslash \GL(V_0)\times \GL(V_1)$ is
generated by $\tr(X)$ and $\langle j X i, \chi\rangle$, where $X$ are
compositions of $B_1$, $B_2$ and $d$ in various orders (so that
$X\colon V_\alpha\to V_\alpha$ in the first case, $X\colon V_0\to V_1$
in the second case), and $\chi\in W^*$. The same statement holds for
$M_0(r,n)$ \cite{Lu-qv}. Therefore the morphism is a closed embedding.
\begin{NB}
  Also two definition gives the same morphism. 
\end{NB}
We have the inverse map given by
\(
   [(B'_1,B'_2,i',j')] \mapsto
   [(B'_1,B'_2,d=\operatorname{id},i',j')]
\)
composed with the extension to $V_0$ or $V_1$ by $0$. Therefore
$\mu^{-1}(0)\dslash \GL(V_0)\times \GL(V_1) \cong M_0(r,n)$.

Composing the natural projective morphism
\(
    \bM_\zeta^{\mathrm{s}}(r,k,n)\to
    \mu^{-1}(0)\dslash \allowbreak\GL(V_0)\times \GL(V_1),
\)
we have a morphism
\begin{equation*}
  \bM_\zeta^{\mathrm{s}}(r,k,n)\to M_0(r,\min(\dim V_0,\dim V_1)).
\end{equation*}

When the parameter $\zeta$ corresponds to the framed moduli space of
torsion free sheaves (${}^\infty\zeta$ in Figure~\ref{fig:zeta1}), the
morphism is one constructed in \cite[Th.3.3]{NY2}.
\begin{NB}
Suppose $0 \le k = \dim V_1 - \dim V_0 \le r$, then
we have $\min(\dim V_0,\dim V_1) = \dim V_0 = 
n - k(r-k)/(2r)$.
\end{NB}%
It factors through the Uhlenbeck (partial) compactification
$\bM_0(r,k,n)$, which corresponds to a parameter on the line
$\zeta_0+\zeta_1 = 0$.

\begin{NB}
\subsection{Examples}

Let us give some examples of ADHM descriptions of torsion free sheaves
on $\bp$.

(1) $V_0 = W = \C$, $V_1 = 0$, $i = 1$. The stability condition
    (S1),(S2) is satisfied, so it gives a line bundle $L$. We have
    $c_1(L) = - [C]$, and hence $L = \shfO(-C)$. In fact, $L$ is given
    by the following display of the monad:
    \begin{equation*}
    \begin{CD}
    @. @. 0 @. 0 
\\
    @. @. @VVV @VVV @.
\\
    0 @>>> O(C-\linf) @>>> \Ker\beta @>>> L @>>> 0
\\
    @. @| @VVV @VVV @.
\\
    0 @>>> O(C-\linf)
    @>{
     \alpha = 
    \lsp{t}{\begin{bmatrix}
     z & w & 0
     \end{bmatrix}}}>>
   \begin{matrix}
    \C^2\otimes\shfO \\ \oplus \\ \shfO
   \end{matrix}
   @>>> 
   \begin{matrix}
   \Coker\alpha 
   \end{matrix}
   @>>> 0
\\
   @. @.
   @V{
   \beta =
   \begin{bmatrix}
    z_2 & - z_1 & i z_0\\
   \end{bmatrix}}VV
   @VVV @.
\\
   @. @. \shfO(\linf) @= \shfO(\linf) @.
\\
   @. @. @VVV @VVV @.
\\
   @. @. 0 @. 0 @.
    \end{CD}
    \end{equation*}

(2) $V_1 = W = \C$, $V_0 = 0$, $j = 1$. The stability condition
    (S1),(S2) is satisfied, so it gives a line bundle $L$. We have
    $c_1(L) = [C]$, and hence $L = \shfO(C)$. It is the dual bundle of
    Example~(1). The display is
    \begin{equation*}
    \begin{CD}
    @. @. 0 @. 0 
\\
    @. @. @VVV @VVV @.
\\
    0 @>>> O(-\linf) @>>> \Ker\beta @>>> L @>>> 0
\\
    @. @| @VVV @VVV @.
\\
    0 @>>> O(-\linf)
    @>{
     \alpha = 
    \lsp{t}{\begin{bmatrix}
     z_1 & z_2 & z_0
     \end{bmatrix}}}>>
   \begin{matrix}
    \C^2\otimes\shfO \\ \oplus \\ \shfO
   \end{matrix}
   @>>> 
   \begin{matrix}
   \Coker\alpha 
   \end{matrix}
   @>>> 0
\\
   @. @.
   @V{
   \beta =
   \begin{bmatrix}
    w & - z & 0\\
   \end{bmatrix}}VV
   @VVV @.
\\
   @. @. \shfO(-C+\linf) @= \shfO(-C+\linf) @.
\\
   @. @. @VVV @VVV @.
\\
   @. @. 0 @. 0 @.
    \end{CD}
    \end{equation*}

(3) The line bundle $\shfO(kC)$ is more complicated. We have
    $\dim V_0 = \frac{k(k-1)}2$, $\dim V_1 = \frac{k(k+1)}2$. For
    example, $k = -3$, we have
\begin{equation*}
\xymatrix @R=.5pc {
W \ar[rd]^i & \\
& V_0(0,0) \ar[rd]^d & \\
& & V_1(0,0) \ar[r]^{B_1} \ar[dd]_{B_2} & V_0(1,0) \ar[rd]^d \\
& & & & V_1(1,0) \ar[r]^{B_1} \ar[dd]^{B_2} & V_0(2,0) \\
& & V_0(0,1) \ar[rd]^d \\
& & & V_1(0,1) \ar[r]^{B_1} \ar[dd]_{B_2}& V_0(1,1) \\
\\
& & & V_0(0,2)
}
\end{equation*}
where $V_\alpha(p,q)$ means that it has weight $t_1^p t_2^q$ if we
consider the data as a $T^2$-fixed point (see below). All vector
spaces are $1$-dimensional, and all the displayed maps are nonzero. We
cannot set either of $V_0(2,0)$, $V_0(1,1)$, $V_0(0,2)$ to be $0$, as
$S_0 = V_0(2,0)\oplus V_0(1,1)\oplus V_0(0,2)$ and $S_1 = V_1(1,0)
\oplus V_1(0,1)$ must satisfy $\dim S_0 > \dim S_1$. Thus the
shape of the diagram is automatically determined by $k$, thanks to the
stability condition. This is different from the example below.

If one of $V_0(2,0)$, $V_0(1,1)$, $V_0(0,2)$ is $0$, the data
satisfies (S1)'. Thus gives a torsion free sheaf. (What are they ?)

(4) Let us try to find the ideal sheaf $I_{p_1}$ at the fixed point
$p_1 = ([1:0:0],[1:0])$. It has $c_1 = 0$, $c_2 = 1$. Therefore $\dim
V_0 = \dim V_1 = 1$. It is a fixed point, so $V_0$, $V_1$ are
$T^2$-modules. If $i = 0$, then $(T_0, T_1) = (0,0)$ violates the
condition~(S2). If $j\neq 0$, then $V_1$ has weight $(-1,-1)$. Then we
must have $B_1 = B_2 = d = 0$ as the weights do not match. But the
equation $B_1 d B_2 - B_2 d B_1 + ij = 0$ is not satisfied. Therefore
we must have $j = 0$. If $B_1 = B_2 = 0$, then $(S_0,S_1) = (0,V_1)$
violates the condition~(S1)'. Therefore either of $B_1$ or $B_2$ is
nonzero. If one is nonzero, the other is zero by the weight condition.
By these discussion, we have two possibilities either $B_1 = 0$ or
$B_2 = 0$.  These correspond to $p_1$ and another fixed point $p_2$.
\begin{equation*}
\xymatrix{
W \ar[rd]^i & \\
V_1(-1,0) \ar[r]_{B_1} & V_0(0,0)
}
\qquad\qquad
\xymatrix{
W \ar[rd]_i & V_1(0,-1) \ar[d]^{B_2} \\
& V_0(0,0)
}
\end{equation*}
Consider the corresponding complex for the first case:
\begin{equation*}
\begin{CD}
\begin{matrix}
   \shfO(C-\linf)
\\
   \oplus
\\ \shfO(-\linf)
\end{matrix}
@>{\alpha}>>
\begin{matrix}
   \C^2\otimes \shfO \\ \oplus \\ \C^2\otimes\shfO \\
   \oplus \\ \shfO
\end{matrix}
@>{\beta}>>
\begin{matrix}
  \shfO(\linf)
\\
  \oplus
\\
  \shfO(-C+\linf)
\end{matrix},
\end{CD}
\end{equation*}
with
\begin{equation*}
   \alpha = 
   \begin{bmatrix}
     z & z_0 \\
     w & 0 \\
     0 & z_1 \\
     0 & z_2 \\
     0 & 0
   \end{bmatrix},
\qquad
   \beta =
   \begin{bmatrix}
    z_2 & - z_1 & 0 & - z_0 & z_0\\
    0   & 0     & w & -z    & 0
   \end{bmatrix}.
\end{equation*}
The linear map induced by $\alpha$ has kernel at $([z_0:z_1:z_2],
[z:w]) = ([1:0:0], [1:0]) = p_1$. Therefore this corresponds to
$I_{p_1}$. The other data gives to $I_{p_2}$.

(5) The ideal sheaf $I_{p_1 + p_2}$ is a `combination' of the two data
in the above example:
\begin{equation*}
\xymatrix{
  W \ar[rd]^i & V_1(0,-1) \ar[d]_{B_2}
\\
  V_1(-1,0) \ar[r]_{B_1} & V_0(0,0),
}
\end{equation*}
where $\dim V_0(0,0) = 2$, $\dim V_1(0,-1) = \dim V_1(-1,0) = 1$. If
the images of $B_1$ and $B_2$ coincide, then $S_0 = \mathrm{image}$,
$S_1 = V_1$ violates the condition~(S1)'. Therefore we may normalize
$B_1 = 
\left[\begin{smallmatrix}
  1 \\ 0
\end{smallmatrix}
\right],
$
$B_2 = 
\left[\begin{smallmatrix}
  0 \\ 1
\end{smallmatrix}
\right].
$
It is also easy to see that the image of $i$ is not either
$
\left[\begin{smallmatrix}
  * \\ 0
\end{smallmatrix}
\right]
$
or
$
\left[\begin{smallmatrix}
  0 \\ *
\end{smallmatrix}
\right].
$
Therefore we may normalize as
\(
   i =
\left[\begin{smallmatrix}
  1 \\ 1
\end{smallmatrix}
\right].
\)

(6) The $T^2$-invariant ideal sheaf $I_Z$ supported at $p_1$
corresponds to the ideal $(x^2,y)$ in the coordinate $(x,y) =
(z_1/z_0, w/z)$. We have $H^0(\shfO/I_Z) = \C 1 \oplus \C x
= \C 1 \oplus \C z_1/z_0$. The weights are $(0,0)$, $(1,0)$
respectively. Therefore we have
\begin{equation*}
\xymatrix @R=0pt {
   W \ar[rdd]^i &
\\ &
\\
   V_1(-1,0) \ar[r]_{B_1} & V_0(0,0) \ar[rdd]^d \\
   || & || \\
   \C\cdot 1/z & \C\cdot 1 & V_1(0,0) \ar[r]^{B_1} & V_0(1,0) \\ 
   & & || & || \\ 
& & \C\cdot z_1/z & \C\cdot x.
}
\end{equation*}
On the other hand, the ideal sheaf corresponding to $(x,y^2)$ is given
by
\begin{equation*}
\xymatrix{
& W \ar[rd]^{i} \\
& V_1(-1,0) \ar[r]_{B_1} \ar[d]^{B_2} & V_0(0,0) \\
V_1(-2,1) \ar[r]^{B_1} & V_0(-1,1)
}
\end{equation*}

The ideal sheaf corresponding to $(x^2, y^2)$ is
\begin{equation*}
\xymatrix {
& W \ar[rd]^i \\
& V_1(-1,0) \ar[r]^{B_1} \ar[d]_{B_2} & V_0(0,0) \ar[rd]^{d} \\
V_1(-2,1) \ar[r]^{B_1} & V_0(-1,1) \ar[rd]^{d} && V_1(0,0)
\ar[r]^{B_1} \ar[d]^{B_2}
& V_0(1,0) \\
&& V_1(-1,1) \ar[r]^{B_1} & V_0(0,1) \\
}
\end{equation*}

A $T^2$-invariant ideal sheaf supported at $p_1$ corresponds to a
Young diagram \cite[2.2]{part1}. The generalization of the above
examples is given by placing $x^p y^q$ on $V_0(p-q,q)$, and also add
$V_1(p-q-1,q)$.
\end{NB}

\section{Moduli spaces}

In this paper we assume $\zeta = (\zeta_0,\zeta_1)\in\R^2$ satisfies
\begin{equation}\label{eq:ass}
  \zeta_0 + \zeta_1 < 0, \quad \zeta_0 < 0
\end{equation}
except in \subsecref{subsec:bp}


\subsection{Smoothness of the moduli space}

We show that $\bM_\zeta^{\mathrm{s}}(r,k,n)$ is smooth in this
subsection.

\begin{Lemma}\label{lem:S2}
  Under the assumption \eqref{eq:ass}, $\zeta$-semistability implies
  the condition~\textup{(S2)}.
\end{Lemma}

\begin{proof}
As $\zeta_1 < -\zeta_0$, we have
\begin{equation*}
     0 \le \zeta_0 \codim T_0 + \zeta_1 \codim T_1 \le
     \zeta_0(\codim T_0 - \codim T_1).
\end{equation*}
As $\zeta_0 < 0$, this implies $\codim T_0 \le \codim T_1$. If the
equality holds, we have $\codim T_1 = 0$ from the second
inequality. Therefore $(T_0,T_1) = (V_0,V_1)$.
\end{proof}

\begin{Lemma}\label{lem:surj}
  Suppose that $(B_1,B_2,d,i,j)$ satisfies \textup{(S2)}. Then $V_0 =
  \Ima B_1 + \Ima B_2 + \Ima i$. The same holds if $(B_1,B_2,d,i,j)$
  is $\zeta$-semistable thanks to \lemref{lem:S2}.
\end{Lemma}

\begin{proof}
Take $T_0 = \Ima B_1 + \Ima B_2 + \Ima i$, $T_1 = V_1$. Then they
satisfy the assumption in (S2). As $\codim T_1 = 0$, we must have
$\codim T_0 = 0$, i.e., $T_0 = V_0$.
\end{proof}

Recall $\mu(B_1,B_2,d,i,j) =  B_1 d B_2 - B_2 d B_1 + ij$. 
\begin{Lemma}\label{lem:stable}
Suppose $(B_1,B_2,d,i,j)$ is $\zeta$-stable.

\textup{(1)} The stabilizer of $(B_1,B_2,d, i, j)$ in $\GL(V_0)\times
\GL(V_1)$ is trivial.

\textup{(2)} The differential $d\mu$ of $\mu$ at $(B_1,B_2,d,i,j)$ is
surjective.
\end{Lemma}

\begin{proof}
(1) Suppose that the pair $(g_0, g_1)\in \GL(V_0)\times
\GL(V_1)$ stabilizes $(B_1,B_2,d,i,j)$. Then $S_0 = \Ima (g_0 -
\id_{V_0})$, $S_1 = \Ima(g_1 - \id_{V_1})$ satisfies the assumption in
(1) of the $\zeta$-stability condition. Therefore we have
\(
    \zeta_0 \dim S_0 + \zeta_1 \dim S_1 \le 0.
\)
Similarly $T_0 = \Ker (g_0 -
\id_{V_0})$, $T_1 = \Ker(g_1 - \id_{V_1})$ satisfy
\(
  \zeta_0 \codim T_0 + \zeta_1 \codim T_1 \ge 0.
\)
But we have $\codim T_0 = \dim S_0$, $\codim T_1 = \dim S_1$, so both
the inequalities must be equalities. Therefore we have $S_0 = 0$, $S_1
= 0$. It means that $g_0 = \id_{V_0}$, $g_1 = \id_{V_1}$.

(2) Suppose $C\in \Hom(V_0,V_1)$ is orthogonal to the image of $d\mu$,
with respect to the pairing given by the trace. Then we have
\begin{equation*}
  CB_\alpha d = d B_\alpha C \ (\alpha=1,2), \quad
  B_1 C B_2 = B_2 C B_1, \quad
  jC = 0, \quad Ci = 0.
\end{equation*}
Let us consider $S_0 = \Ima(B_1 C)$, $S_1 = \Ima (C B_1)$,
$T_0 = \Ker(B_1 C)$, $T_1 = \Ker (C B_1)$. By the same argument as in (1),
the $\zeta$-stability implies $S_0 = S_1 = 0$, that is
\begin{equation*}
   B_1 C = 0, \quad C B_1 = 0.
\end{equation*}
Exchanging $B_1$ and $B_2$, we also have 
\(
 B_2 C = 0, 
\)
\(
C B_2 = 0.
\)
Then $C$ is equal to $0$ on $\Ima B_1 + \Ima B_2 + \Ima i$. Then
by \lemref{lem:surj} we have $C = 0$.
\end{proof}

By the standard argument from the geometric invariant theory, we have
\begin{Theorem}
$\bM_\zeta^{\mathrm{s}}(r,k,n)$ is nonsingular of dimension
\begin{equation*}
   \dim W (\dim V_0 + \dim V_1) - (\dim V_0 - \dim V_1)^2 
   = 2nr,
\end{equation*}
provided $\bM_\zeta^{\mathrm{s}}(r,k,n)\neq\emptyset$. Moreover,
$\bM_\zeta^{\mathrm{s}}(r,k,n)$ is the fine moduli space of
$\zeta$-stable $\GL(V_0)\times \GL(V_1)$-orbits.
\end{Theorem}

For a later purpose (in the subsequent paper) let us describe the
tangent space: it is the middle cohomology of
\begin{equation}\label{eq:tangentcpx}
\begin{CD}
\begin{matrix}
   \Hom(V_0,V_0)
\\
   \oplus
\\ \Hom(V_1,V_1)
\end{matrix}
@>{\iota}>>
\begin{matrix}
   \Hom(V_0,V_1)\\ \oplus \\ \C^2\otimes \Hom(V_1, V_0) \\ 
   \oplus \\ \Hom(W,V_0) \\ \oplus \\ \Hom(V_1,W)
\end{matrix}
@>{d\mu}>>
\begin{matrix}
  \Hom(V_1,V_0)
\end{matrix},
\end{CD}
\end{equation}
with
\begin{equation*}
\begin{split}
  \iota
   \begin{bmatrix}
     \xi_0 \\ \xi_1
   \end{bmatrix}
   = 
   \begin{bmatrix}
     d \xi_0 - \xi_1 d     \\
     B_1 \xi_1 - \xi_0 B_1 \\
     B_2 \xi_1 - \xi_0 B_2 \\
     \xi_0 i \\
      - j \xi_1
   \end{bmatrix},
   (d\mu)
   \begin{bmatrix}
    \widetilde d \\ \widetilde B_1 \\ \widetilde B_2 \\ \widetilde i \\
    \widetilde j
   \end{bmatrix}
   = 
   \begin{aligned}
   & B_1 d \widetilde B_2 + B_1 \widetilde d B_2 + \widetilde B_1 d
   B_2
\\
   & \quad
   - B_2 d \widetilde B_1 - B_2 \widetilde d B_1 - \widetilde B_2 d
   B_1
\\
   & \quad\quad
   + \widetilde i j + i\widetilde j,
   \end{aligned}
\end{split}
\end{equation*}
where $d\mu$ is the differential of $\mu$, and $\iota$ is the
differential of the group action. Remark that $d\mu$ is
surjective and $\iota$ is injective by the above lemma.

\subsection{The $\zeta$-stable point with $W=0$}

For a later purpose we study the case $W=0$ in this subsection. We
assume \eqref{eq:ass} and
\begin{equation}\label{eq:assW=0}
  \zeta_0 \dim V_0 + \zeta_1 \dim V_1 = 0.
\end{equation}

\begin{Definition}\label{def:stableW=0}
  Suppose $W=0$ and $\zeta_0 \dim V_0 + \zeta_1 \dim V_1 = 0$. We say
  $X = (B_1, B_2, d)$ as above is {\it $\zeta$-semistable\/} if the
  condition (1) in Definition~\ref{def:stable} holds. (Note that the
  condition (2) is equivalent to (1) when $W=0$). And the {\it
    $\zeta$-stability\/} is defined by requiring that the inequality
  is strict unless $(S_0,S_1) = (0,0)$ or $(V_0,V_1)$.
\end{Definition}

We consider the quotient space
\begin{equation*}
    \left\{ (B_1,B_2,d) \left|
        \begin{aligned}[m]
    & B_1 d B_2 - B_2 d B_1 = 0 \\
    & \text{$\zeta$-stable}      
        \end{aligned}
  \right\}\right/ \left(\GL(V_0)\times \GL(V_1)/\C^*\right),
\end{equation*}
where $\C^* = \{ (\lambda\id_{V_0},\lambda \id_{V_1}) \in
\GL(V_0)\times \GL(V_1) \mid \lambda\in\C^*\}$. As it acts trivially
on the data, the quotient group acts on the space.

\begin{Lemma}\label{lem:d=0}
Suppose $(B_1,B_2,d)$ is $\zeta$-stable. Then $d = 0$.  
\end{Lemma}

\begin{proof}
Consider $S_0 = \Ker(B_1 d)$, $S_1 = \Ker(d B_1)$. By the equation
$B_1 d B_2 = B_2 d B_1$, they satisfy the assumption in the
$\zeta$-semistability condition, hence we have $\zeta_0 \dim S_0 +
\zeta_1 \dim S_1 \le 0$. Similarly we have $\zeta_0 \dim S'_0 +
\zeta_1 \dim S'_1 \le 0$ for $S'_0 = \Ima(B_1 d)$, $S'_1 = \Ima (d
B_1)$. But as $\dim S_\alpha + \dim S'_\alpha = \dim V_\alpha$
($\alpha=0,1$), the assumption \eqref{eq:assW=0} implies that both
inequalities must be equalities. Hence $S_0 = 0$, $S_1 = 0$ or $S_0 =
V_0$, $S_1 = V_1$. Suppose that the first case occurs. Then both
$B_1d$ and $d B_1$ are injective. Therefore both $B_1$, $d$ must be
isomorphisms, and hence $\dim V_0 = \dim V_1$. But this is impossible
by \eqref{eq:ass} and \eqref{eq:assW=0}. Therefore $B_1 d = 0$, $d B_1
= 0$. Exchanging $B_1$ and $B_2$, we have $B_2 d = 0$, $d B_2 = 0$.
But as we have $V_0 = \Ima B_1 + \Ima B_2$ by the same argument as in
\lemref{lem:surj}, we get $d=0$.
\end{proof}

Therefore the data $(B_1,B_2,d) = (B_1,B_2,0)$ can be considered as a
representation of the Kronecker quiver. By the same argument in
\lemref{lem:stable}(1), it has the trivial stabilizer in
$\GL(V_0)\times\GL(V_1)/\C^*$. In particular, it is an indecomposable
representation. The classification of indecomposable representations
of the Kronecker quiver is well-known (see e.g., \cite[1.8]{GR}). If we
take suitable bases of $V_0$ and $V_1$, the data is written as either
of the followings:
\begin{subequations}\label{eq:Kronecker}
\begin{align}
  & B_1 = \left[\begin{smallmatrix} 1_m \\ 0
    \end{smallmatrix}\right],\quad
  B_2 = \left[\begin{smallmatrix} 0 \\ 1_m
    \end{smallmatrix}\right]\in \operatorname{Mat}(m+1,m)\quad
  \text{($m\ge 0$)},
\\
  & B_1 = \left[\begin{smallmatrix} 1_m & 0
    \end{smallmatrix}\right], \quad
  B_2 = \left[\begin{smallmatrix} 0 & 1_m
    \end{smallmatrix}\right]\in \operatorname{Mat}(m,m+1)\quad
  \text{($m\ge 0$)},
\\
  & B_1 = 1_m, \quad B_2 = J_m \quad\text{($m\ge 1$)},
\\
  & B_1 = a1_m + J_m, \quad B_2 = 1_m \quad\text{($a\in\C$, $m\ge 1$)}.
\end{align}
\end{subequations}
Here $1_m$ is the identity matrix of size $m$, $J_m$ is the Jordan
block of size $m$ with eigenvalue $0$, and $\operatorname{Mat}(k,l)$
is the space of all $k\times l$-matrices.
By \eqref{eq:ass} and \eqref{eq:assW=0} we have $\dim V_0 < \dim V_1$.
Thus we only have the possibility b).

\begin{Remark}
It is well-known that the derived category of finite dimensional
representations of the Kronecker quiver is equivalent to the derived
category of coherent sheaves on $\proj^1$. The equivalence is given by
\begin{equation*}
    R\Hom_{\proj^1}(\shfO_{\proj^1}\oplus \shfO_{\proj^1}(1),\bullet)\colon
    D^b(\mathrm{Coh}\ \proj^1) \to D^b(\textrm{mod}\mathchar`-A),
\end{equation*}
where $A$ is the path algebra of the Kronecker quiver, or more
explicitly, $V_0 = R\Hom_{\proj^1}(\shfO_{\proj^1},\bullet)$, $V_1 =
R\Hom_{\proj^1}(\shfO_{\proj^1}(1),\bullet)$ with homomorphisms given by
$z_1$ and $z_2$ of the homogeneous coordinate $[z_1:z_2]$ of
$\proj^1$. Then the objects corresponding to above indecomposable
representations are
\begin{aenume}
\item $\shfO(m)$ ($m\ge 0$),
\item $\shfO(-m-1)[1]$ ($m\ge 0$),
\item $\shfO/z_2^m \shfO$ ($m\ge 1$),
\item $\shfO/(z-a)^m\shfO$ ($z=z_1/z_2$, $a\in\C$, $m\ge 1$).
\end{aenume}%
\begin{NB}
In particular, the underlying abelian category is different, and
$\shfO(-m-1)$ ($m\ge 0$) is `moved'.
\end{NB}
\end{Remark}

\begin{Lemma}\label{lem:C_mstable}
The data given in \textup{b)} is $\zeta$-stable.
\end{Lemma}

\begin{proof}
We take $V_0 = \C^m$, $V_1 = \C^{m+1}$ and $B_1$, $B_2$ as in b). We
have
$m \zeta_0 + (m+1)\zeta_1 = 0$ by the assumption \eqref{eq:assW=0}. We
further have $\zeta_1\ge 0$ by \eqref{eq:ass}.

The condition $S_0\subset V_0$, $S_1\subset V_1$ with
$B_\alpha(S_1)\subset S_0$ means that the restriction
$(\left.B_1\right|_{S_0}, \left.B_2\right|_{S_0})$ to $S_0$, $S_1$ is 
a subrepresentation of $(B_1,B_2)$ as the representation of the
Kronecker quiver.
By the well-known representation theory of the Kronecker quiver, or by
a direct calculation, we know that a quotient representation of $B_1$,
$B_2$ is a direct sum of indecomposable representations $(B_1',
B_2')$, $(B_1'',B_2'')$,\dots, all of type b) with size $n'$, $n''$,\dots.
with $n' + n'' + \cdots \le m$. Then we have
\begin{equation*}
\begin{split}
   & \zeta_0 \codim S_0 + \zeta_1 \codim S_1
\\
  = \; & \zeta_0(n' + n'' + \cdots) + \zeta_1 ((n' + 1) + (n'' + 1) + \cdots)
\\
  = \; &
    (\zeta_0 + \zeta_1)\left\{(n'+n''+\cdots) - m\right\} 
    + \zeta_1 (\#\{ \mathrm{factors}\} - 1) \ge 0,
\end{split}
\end{equation*}
where we have used $m\zeta_0 + (m+1)\zeta_1 = 0$ in the second
equality and $\zeta_0 + \zeta_1 < 0$ and $\zeta_1 \ge 0$ in the last
inequality.
It is also clear that we have $(S_0, S_1) = (0,0)$ or $(V_0,V_1)$ if
the equality holds.%
\begin{NB}
  We implicitly assume $(S_0,S_1) \neq (V_0,V_1)$ and hence $\#\{
  \mathrm{factors}\} \ge 1$, as the assertion becomes trivial
  otherwise. Suppose $\zeta_1 = 0$. Then $m=0$ by the assumption
  $m\zeta_0 + (m+1)\zeta_1 = 0$. The assertion is trivial in this
  case.  So we may assume $\zeta_1\neq 0$. Then the equality holds if
  and only if $(n'+n''+\cdots) = m$ and $\#\{ \mathrm{factors}\} = 1$.
  This means $(S_0,S_1) = (0,0)$.
\end{NB}
\end{proof}

In summary we have
\begin{Theorem}\label{thm:W=0}
Assume the condition \eqref{eq:ass}.
The moduli space of $\zeta$-stable data $X = (B_1,B_2,d)$ is a single
point given by the data \textup{b)} above with $d=0$ if $\dim V_0 =
m$, $\dim V_1 = m+1$ with $m\zeta_0 + (m+1)\zeta_1 = 0$.  And it is
the empty set otherwise.
\end{Theorem}

Let $C_m = (B_1,B_2,0)$ denote the data given in this theorem.
In \propref{prop:Cm} we prove that $C_m$ corresponds to
$\shfO_C(-m-1)$.

\begin{NB}
The following lemma was originally used to classify the $\zeta$-stable
data. But it is, in fact, unnecessary as we have done above. 

\begin{Lemma}\label{lem:dimension}
Suppose $(B_1,B_2,d)$ is $\zeta$-stable. Then

\textup{(1)} The stabilizer of $(B_1,B_2,d)$ in $\GL(V_0)\times
\GL(V_1)/\C^*$ is trivial.

\textup{(2)} The differential of the map
$(B_1,B_2,d)\mapsto B_1 d B_2 - B_2 d B_1$ at $(B_1,B_2,d)$ is surjective

\textup{(3)} $\dim V_0 = \dim V_1 - 1$.  
\end{Lemma}

\begin{proof}
By the similar argument as in
the proof of \lemref{lem:stable}(1), the statement (1) can be proved.

Next suppose that $C$ is orthogonal to the image of the differential
of the map $(B_1,B_2,d)\mapsto B_1 d B_2 - B_2 d B_1$. Then we have
\begin{equation*}
  CB_\alpha d = d B_\alpha C \ (\alpha=1,2), \quad
  B_1 C B_2 = B_2 C B_1.
\end{equation*}
Let us consider $S_0 = \Ima(B_1 C)$, $S_1 = \Ima (C B_1)$, $S_0' =
\Ker(B_1C)$, $S_1' = \Ker(C B_1)$.  By the same argument as in
\lemref{lem:stable}(2), the $\zeta$-stability implies $S_0 = S_1 = 0$
or $S_0 = V_0$, $S_1 = V_1$. Suppose we have the latter case. Then we
have both $B_1C$ and $C B_1$ must be surjective, hence both $B_1$ and
$C$ are isomorphisms. In particular, we must have $\dim V_0 = \dim
V_1$. But then we have a contradiction between \eqref{eq:ass} and
\eqref{eq:assW=0}. So we must have $S_0 = 0$, $S_1 = 0$, i.e.,
$B_1 C = 0$, $C B_1 = 0$. Similarly we have $B_2 C = 0$, $C B_2 =
0$. But as in \lemref{lem:surj} we have
$V_0 = \Ima B_1 + \Ima B_2$. Therefore we must have $C = 0$.
Hence the quotient space is smooth of dimension
\begin{equation*}
  1 - (\dim V_0 - \dim V_1)^2
\end{equation*}
provided it is nonempty. In particular, we must have
\(
  \dim V_0 - \dim V_1 = 0 
\)
or 
\(
   \pm 1.
\)
However $\dim V_0 - \dim V_1$ cannot be $0$ or $1$ from the
assumptions \eqref{eq:ass} and \eqref{eq:assW=0}. Hence it must be
equal to $-1$.
\end{proof}
\end{NB}

\subsection{Blowup of the plane as a moduli space}\label{subsec:bp}

If the condition \eqref{eq:ass} is not satisfied, \thmref{thm:W=0} is
no longer true. We consider the limiting case 
\begin{equation}
  \label{eq:ass'}
 \zeta_0 + \zeta_1 = 0, \qquad \zeta_0 < 0.  
\end{equation}
From \eqref{eq:assW=0} we have $\dim V_0 = \dim V_1$. The data
$(B_1,B_2,d)$ is $\zeta$-stable if and only if it satisfies
\begin{description}
\item[(S0)] for subspaces $S_0\subset V_0$, $S_1\subset V_1$ such that
  $B_\alpha(S_1)\subset S_0$ ($\alpha=1,2$), $d(S_0)\subset S_1$, we
  have either $\dim S_0 > \dim S_1$, $S_0 = S_1 = 0$, or $S_0 = V_0$,
  $S_1 = V_1$.
\end{description}
This is the `$W=0$' version of the conditions (S1), (S2) in
\thmref{thm:King}.

\begin{Lemma}
Suppose $B_1$, $B_2$, $d$ satisfies \textup{(S0)} and
$B_1 d B_2 = B_2 d B_1$.
Then $\dim V_0 = \dim V_1 = 1$.
\end{Lemma}

\begin{proof}
From the equation $B_1 d B_2 = B_2 d B_1$, $B_1 d$ and $B_2 d$
commute. Let $0\neq v_0\in V_0$ be a simultaneous eigenvector, and let
$v_1 = d v_0$. Then $(S_0, S_1) = (\C v_0, \C v_1)$ satisfies the
assumption in (S0). As we cannot have $(S_0,S_1) = (0,0)$, we either
have $(S_0,S_1) = (V_0,V_1)$ or $v_1 = 0$. We are done in the first
case. Therefore we may assume the latter case. Hence $\Ker d\neq 0$.

We take $S_0$, $S_1$, $S'_0$, $S'_1$ as in the proof of
\lemref{lem:d=0}. By the same argument there, we get
$(S_0,S_1) = (0,0)$ or $(V_0,V_1)$. But the first case cannot occur
as $\Ker d\neq 0$. Therefore only the latter case can occur, thus
$d B_1 = 0$. Exchanging $B_1$ and $B_2$, we get $d B_2 = 0$. By the
same argument as in \lemref{lem:surj}, the condition (S0) implies
$V_0 = \Ima B_1 + \Ima B_2$. Therefore we have $d=0$. Therefore
$(B_1,B_2,d) = (B_1,B_2,0)$ can be regarded as a representation of the
Kronecker quiver. Moreover, by the same argument as in
\lemref{lem:stable}(1), it is indecomposable. Therefore after taking
a base, $(B_1,B_2)$ can be written as in \eqref{eq:Kronecker}. As
$\dim V_0 = \dim V_1$, the cases a),b) cannot happen. In the case c)
or d), we take
\begin{equation*}
   v_1 =
   \left[\begin{smallmatrix}
     1 \\ 0 \\ \vdots \\ 0
   \end{smallmatrix}\right].
\end{equation*}
Then $B_1 v_1$ and $B_2 v_1$ are linearly dependent, and span a
$1$-dimensional subspace $S''_0$ in $V_0$. As the pair $(S''_0,S''_1 =
\C v_1)$ violates the inequality in (S0), we must have $S''_0 = V_0$,
$S''_1 = V_1$. Thus $\dim V_0 = \dim V_1 = 1$ also in this case.
\end{proof}

\begin{Theorem}\label{thm:blowupplane}
Assume \eqref{eq:ass'} and $\dim V_0 = \dim V_1$.
Then the moduli space of $\zeta$-stable data $(B_1,B_2,d)$ is
isomorphic to the blowup of the plane at the origin $\widehat{\mathbb
  C}^2$ if $\dim V_0 = \dim V_1 = 1$.  And it is the empty set
otherwise.
\end{Theorem}

\begin{proof}
By the previous lemma, we may assume $V_0 = V_1 = \C$. Then we have a
morphism from the moduli space of $\zeta$-stable data $(B_1,B_2,d)$ to
$\widehat{\mathbb C}^2$ given by
\begin{equation*}
\begin{split}
    & [(B_1,B_2,d)] \mapsto
    ((B_1d, B_2 d), [B_1:B_2])
    \in 
    \widehat{\mathbb C}^2;
\\
  & \qquad\text{where }
 \widehat{\mathbb C}^2 =
    \left\{ ((z_1,z_2),[z:w])\in\C^2\times\proj^1 \mid z_1 w = z_2 z \right\}.
\end{split}
\end{equation*}
As we observed already, we have $V_0 = \Ima B_1 + \Ima B_2$ if
$(B_1,B_2,d)$ satisfies (S0) by the argument in \lemref{lem:surj}.
Therefore $B_1 = B_2 = 0$ is not possible. Therefore
$[B_1:B_2]\in\proj^1$ is defined, and the above actually defines a
morphism from the moduli space to $\widehat{\mathbb C}^2$.

Conversely take a point $((z_1,z_2),[z:w])\in \widehat{\mathbb
  C}^2$. We set $V_0 = V_1 = \C$ and
\begin{equation*}
  B_1 = z, \quad B_2 = w, \quad d =
  s = 
  \begin{cases}
    z_1/z & \text{if $z\neq 0$},
    \\
    z_2/w & \text{if $w\neq 0$}.
  \end{cases}
\end{equation*}
In this case, the condition (S0) is equivalent to $V_0 = \Ima B_1 +
\Ima B_2$, and hence is satisfied. Moreover, the above is well-defined
up to the action of $\GL(V_0)\times\GL(V_1)/\C^* \cong \C^*$. It is
also clear that it is the inverse of the previous map. Hence we have
the assertion.
\end{proof}

From this theorem, we also see that the moduli space
$\bM_\zeta^{\mathrm{ss}}(0,0,n)$ of $S$-equivalence classes of
$\zeta$-semistable data is the $n^{\mathrm{th}}$ symmetric product of
$\widehat{\mathbb C}^2$.

It is instructive to consider what happens if $\zeta_0$ becomes
positive. (We still assume $\zeta_0+\zeta_1 = 0$ and $\dim V_0 = \dim
V_1 = 1$.) In this case, the $\zeta$-stability is equivalent to $\Ker
d = 0$, i.e., $d$ is an isomorphism. Then the moduli space is
isomorphic to $\C^2$ under the map
\begin{equation*}
    [(B_1,B_2,d)] \mapsto
    (B_1d, B_2 d) \in \C^2.
\end{equation*}
The converse is given by 
\[
   (z_1,z_2) \mapsto (B_1,B_2,d) = (z_1,z_2,1).
\]
When we cross the wall $\zeta_0 = 0$, the exceptional locus $C = \{
[(B_1,B_2,0)] \mid [B_1:B_2]\in\proj^1\}$ is replaced by the single
point $0 = \{ [(0,0,1)] \}$. This example is a prototype of the
picture given in the subsequent paper. 

\begin{Remark}\label{rem:allW=0}
  It is straightforward to generalize the above arguments to classify
  all $\zeta$-stable solution with $W=0$ for any parameter $\zeta$
  which satisfies \eqref{eq:assW=0}, but not necessarily
  \eqref{eq:ass}. We leave the proof as an exercise for a reader.

(1) Case $\zeta_0 + \zeta_1 < 0$, $\zeta_0\ge 0$:
The only possibility is $V_0 = \C$, $V_1 = 0$ with $\zeta_0 = 0$.

(2) Case $\zeta_0 + \zeta_1 > 0$, $\zeta_1 > 0$:
The data given by (\ref{eq:Kronecker}a) with $d=0$
if $(m+1)\zeta_0 + m\zeta_1 = 0$.

(3) Case $\zeta_0 + \zeta_1 > 0$, $\zeta_1 \le 0$:
The only possibility is $V_1 = \C$, $V_0 = 0$ with $\zeta_1 = 0$.

(4) Case $\zeta_0 + \zeta_1 = 0$, $\zeta_0 > 0$: We have $\dim V_0 =
\dim V_1 = 1$ and the moduli space is isomorphic to $\C^2$, as
mentioned above.

(5) Case $\zeta_0 = \zeta_1 = 0$: We have $\dim V_0 = \dim V_1 = 1$
and the moduli space is isomorphic to $\C^2\setminus\{(0,0)\}$.

\begin{NB}
  Answer to Exercise: 

  There is a symmetry $(B_1,B_2,d,i,j)\mapsto
  ({}^t B_1,{}^t B_2, {}^t d, {}^tj, {}^t i)$ where the latter is for
  $(V_1^*, V_0^*, W^*)$. (Note that we put $V_1^*$ on the vertex $0$.)
  Under this exchange, a $(\zeta_0,\zeta_1)$-stable object is mapped
  to a $(-\zeta_1,-\zeta_0)$-stable object. Therefore the case (2) is
  reduced to the case we already studied, and the case (3) is
  equivalent to the case (1). In the case (1), \lemref{lem:S2} is
  changed as

  $\zeta$-semistability implies :
  `$T_0\subset V_0$, $T_1\subset V_1$ with $B_\alpha(T_1)\subset
  T_0$, $d(T_0)\subset T_1$ satisfy $\codim T_1 < \codim T_0$ unless
  $(T_0,T_1) = (V_0,V_1)$.'

  Then \lemref{lem:surj} is changed to $\Ima d = V_1$. We follow the
  proof of \lemref{lem:d=0} until we show $B_1 d = 0$, $B_2 d = 0$.
  The surjectivity of $d$ implies $B_1 = B_2 = 0$. Now we use the
  classification of indecomposable representations of the $A_2$-quiver,
  which are equivalent to a linear algebra problem.
\end{NB}
\end{Remark}

\section{Resolution of the diagonal}

\begin{NB}
Jan. 22: I correct the definition of the complex $C^\bullet$. I also
change the definition of $\Qcal$ so that it becomes clear that it is a
pull-back of a vector bundle over $\proj^2$. I also change the last
part of the proof according to Kota's note.
\end{NB}

Let us define a rank $2$ vector bundle $\Qcal$ over $\bp$ as the
quotient
\begin{equation*}
   \Qcal \defeq \shfO^{\oplus 3} / \Ima \iota,
   \qquad\text{where  }
   \iota =
   \begin{bmatrix}
     z_0 \\ z_2 \\ -z_1
   \end{bmatrix}
   \colon \shfO(-\linf) \to \shfO^{\oplus 3}.
\end{equation*}
We have $\Wedge^2\Qcal = \shfO(\linf)$. 
\begin{NB}
  This is the pull-back of the similar bundle over $\proj^2$ used in
  \cite[\S2.1]{Lecture}.
\end{NB}%
This can be also described as the quotient
\begin{equation*}
  0 \to
  \shfO(C-\linf) \xrightarrow{
    \left[
      \begin{smallmatrix}
        z_0 \\ w \\ -z
      \end{smallmatrix}
    \right]
  }
  \shfO(C)\oplus \shfO^{\oplus 2} \to \Qcal \to 0,
\end{equation*}
where the isomorphism is the one induced by
\(
  \operatorname{diag}(s,\operatorname{id},\operatorname{id})
  \colon \shfO \oplus \shfO^{\oplus 2}
  \to \shfO(C) \oplus \shfO^{\oplus 2}.
\)
\begin{NB}
Then
\(
u \left[
  \begin{smallmatrix}
    z_0 \\ z_2\\ -z_1
  \end{smallmatrix}
\right]
\)
is maps to
\(
us \left[
  \begin{smallmatrix}
    z_0 \\ w\\ -z
  \end{smallmatrix}
\right].
\)
Hence 
\(
  \Qcal\longrightarrow 
  \shfO(C)\oplus \shfO^{\oplus 2}
  /\Ima \left(\left[\begin{smallmatrix}
    z_0 \\ w\\ -z
  \end{smallmatrix}\right]\colon
   \shfO(C-\linf)\to \shfO(C)\oplus\shfO^{\oplus 2}\right)
\)
is well-defined. It is an isomorphism over $\bp\setminus \linf$ thanks
to the identification $\mathcal Q|_{\bp\setminus \linf} \cong
\shfO^{\oplus 2}$ below, as $sz = z_1$, $sw = z_2$. And over
$\bp\setminus C$, we trivialize $\shfO(C)$ by setting $s=1$, we see
that the homomorphism is an isomorphism.
\end{NB}
We define a section $\varphi$ of the vector bundle
$\Qcal\boxtimes\shfO(\linf)$ over $\bp\times\bp$ as the composite of
\begin{equation*}
   \varphi\colon \shfO\boxtimes\shfO \xrightarrow{
     \left[\begin{smallmatrix}
       \id\boxtimes z_0\\ \id\boxtimes z_2 \\ -\id\boxtimes z_1
     \end{smallmatrix}\right]}
   \begin{matrix}
    \shfO^{\oplus 3}\boxtimes\shfO(\linf)
   \end{matrix}
    \xrightarrow{p\boxtimes\id_{\shfO(\linf)}}
    \Qcal\boxtimes\shfO(\linf),
\end{equation*}
where the second arrow is the homomorphism induced by the projection
$p\colon \shfO^{\oplus 3}\to\Qcal$. 
%
%
We also define homomorphisms $\psi$, $\chi$ as composites of
\begin{equation*}
  \psi\colon
  \begin{gathered}[t]
   \shfO(-C)\boxtimes\shfO(C)
   \xrightarrow{\left[
       \begin{smallmatrix}
         0 \\ \id\boxtimes w \\ -\id\boxtimes z
       \end{smallmatrix}
     \right]}
   \shfO(-C)^{\oplus 3}\boxtimes\shfO(\linf)
   \to\Qcal(-C)\boxtimes\shfO(\linf),
\\
  \chi\colon
   \shfO\boxtimes\shfO
   \xrightarrow{\left[
       \begin{smallmatrix}
         -\id\boxtimes z_0 \\ 0 \\ 0
       \end{smallmatrix}
     \right]}
   \begin{matrix}
     \shfO \boxtimes\shfO(\linf) \\ \oplus \\
     \shfO(-C)^{\oplus 2} \boxtimes\shfO(\linf)
   \end{matrix}
   \to\Qcal(-C)\boxtimes\shfO(\linf)
  \end{gathered}
\end{equation*}
respectively. Here we used the second description of $\Qcal$ for the
definition of $\chi$.

In order to distinguish the coordinates on the first and second
factors of $\bp$, we denote the former by $([z_0',z_1',z_2'],[z':w'])$
and the latter by $([z_0,z_1,z_2],[z:w])$ hereafter.
We define a complex $C^\bullet$ of vector bundles over
$\bp\times\bp$ by
\begin{multline}\label{eq:diagonal}
C^\bullet\colon
\begin{matrix}
   \shfO\boxtimes \shfO(C-2\linf)
\\
   \oplus
\\ \shfO(C-\linf)\boxtimes\shfO(-2\linf)
\end{matrix}
\\
\xrightarrow{d^{-1}}
\begin{matrix}
   \shfO^{\oplus 2}\boxtimes \shfO(-\linf)
   \\ \oplus \\
   \Qcal^\vee(C)\boxtimes \shfO(-\linf)
\end{matrix}
\xrightarrow{d^0}
\begin{matrix}
  \shfO\boxtimes\shfO \\ \oplus \\
  \shfO(C)\boxtimes\shfO(-C)
\end{matrix}
\end{multline}
with
\begin{equation*}
\begin{split}
   & d^{-1} =
   \begin{bmatrix}
     \id\boxtimes z & z'\boxtimes z_0 \\
     \id\boxtimes w & w'\boxtimes z_0 \\
     0 & \varphi\lrcorner
   \end{bmatrix},
\qquad
d^0 =
   \begin{bmatrix}
    \id\boxtimes z_2 & -\id\boxtimes z_1 & \chi^{\vee}\\
    s'\boxtimes w & - s'\boxtimes z & \psi^{\vee}
   \end{bmatrix}.
\end{split}
\end{equation*}
We assign the degree $0$ to the middle term.
Let $\Delta \equiv \Delta_{\bp}$ be the diagonal in $\bp\times\bp$.
Let $t\colon C^1 \to \shfO_{\Delta}$ be the homomorphism given by
restricting to the diagonal and taking the difference of traces after
the isomorphism
\(
   \left(\shfO\boxtimes\shfO\oplus\shfO(C)\boxtimes\shfO(-C)
   \right)|_{\Delta}
   \cong {\mathcal H}om(\shfO,\shfO)\oplus {\mathcal
   H}om(\shfO(C),\shfO(C)).
\)

Over
\(
 \bp\setminus\linf = \{ ([1:z_1':z_2'],[z':w']) \},
\) 
we can identify $\Qcal$ with $\shfO^{\oplus 2}$ by
\begin{equation*}
\begin{split}
   & \left(\shfO^{\oplus 3}\right)/\shfO(-\linf)
   \ni 
  \begin{bmatrix}
      a \\ b \\ c
  \end{bmatrix}
  \bmod \shfO(-\linf)
  \longmapsto
  \begin{bmatrix}
    b - az_2' \\ c + az_1'
  \end{bmatrix},
\\
  & 
   \left(\shfO(C)\oplus\shfO^{\oplus 2}\right)/\shfO(C-\linf)
   \ni 
  \begin{bmatrix}
      a \\ b \\ c
  \end{bmatrix}
  \bmod \shfO(C-\linf)
  \longmapsto
  \begin{bmatrix}
    b - aw' \\ c + az'
  \end{bmatrix}.
\end{split}
\end{equation*}
The sections $\varphi$, $\psi$, $\chi$ are re-written as
\begin{gather*}
   \varphi\colon\shfO\boxtimes\shfO
   \xrightarrow{\left[
     \begin{smallmatrix}
      \id\boxtimes z_2 - z_2'\boxtimes z_0
      \\
      - \id\boxtimes z_1 + z_1'\boxtimes z_0
     \end{smallmatrix}
     \right]}
     \shfO^{\oplus 2} \boxtimes \shfO(\linf),
\\
   \psi\colon
   \shfO(-C)\boxtimes\shfO(C)
   \xrightarrow{\left[
       \begin{smallmatrix}
         \id\boxtimes w \\ -\id\boxtimes z 
       \end{smallmatrix}
     \right]}
   \shfO(-C)^{\oplus 2} \boxtimes\shfO(\linf),
\\
  \chi\colon
   \shfO\boxtimes\shfO
   \xrightarrow{\left[
       \begin{smallmatrix}
         w'\boxtimes z_0 \\ -z'\boxtimes z_0
       \end{smallmatrix}
     \right]}
   \shfO(-C)^{\oplus 2}\boxtimes\shfO(\linf).
\end{gather*}
\begin{NB}
And
  \begin{equation*}
    \varphi\lrcorner
    \colon \shfO\boxtimes\shfO(-\linf)
    \xrightarrow{\left[
     \begin{smallmatrix}
      \id\boxtimes z_1 - z_1'\boxtimes z_0
      \\
      \id\boxtimes z_2 - z_2'\boxtimes z_0
     \end{smallmatrix}
     \right]}
     \shfO^{\oplus 2} \boxtimes \shfO.
  \end{equation*}
Therefore
\begin{equation*}
  \begin{split}
    & \chi^\vee \circ(\varphi\lrcorner)
    = w'\boxtimes z_0 z_1 - z'\boxtimes z_0 z_2, 
\\
    & \psi^\vee \circ(\varphi\lrcorner)
    = z_2'\boxtimes z_0 z - z_1'\boxtimes z_0 w.
  \end{split}
\end{equation*}
\end{NB}%
We can directly check $d^0\circ d^{-1} = 0$, $t\circ d^0 = 0$ over
\(
 (\bp\setminus\linf) \times\bp.
\)
It holds over $\bp\times\bp$ by the continuity.

\begin{Proposition}[\protect{cf.\ \cite[\S3]{KN}}]\label{prop:diagonal}
The complex
\begin{equation*}
   0\to C^{-1} \xrightarrow{d^{-1}} C^0 
   \xrightarrow{d^0} C^1 \xrightarrow{t} \shfO_{\Delta_\bp}\to 0
\end{equation*}
gives a resolution of $\shfO_{\Delta_{\bp}}$.%
\begin{NB}
More precisely it gives a resolution of $\shfO_{\Delta_\bp}[-1]$.
\end{NB}
\end{Proposition}

In the following lemma we consider homomorphisms $d^0$, $d^{-1}$, etc, as
linear maps between fibers of vector bundles.

\begin{Lemma}\label{lem:diag}
\textup{(1)} $\varphi$ vanishes over $\left(C\times
    C\right)\cup\Delta_{\bp}$.

\textup{(2)} The locus $\{ \text{$d^{-1}$ is not injective}\}$ is
$\Delta_{\bp}$. And $\Ker d^{-1}|_{\Delta_{\bp}} = \{ (-z_0 f)\oplus
f \mid f\in \shfO(C-3\linf)\} \cong \shfO(C-3\linf) \cong K_{\bp}$.

\textup{(3)} $d^0$ is surjective outside $\Delta_{\bp}$, and the
cokernel of $d^0$ is isomorphic to the one-dimensional space of scalar
endomorphisms $\{ \lambda \id_\shfO \oplus\allowbreak (-\lambda)
\id_{\shfO(C)}\}$ on the diagonal $\Delta_{\bp}$ under the isomorphism
\(
   \bigl(\shfO\boxtimes\shfO\oplus\allowbreak
     \shfO(C)\boxtimes\shfO(-C)
   \bigr)|_{\Delta_{\bp}}
   \cong {\mathcal H}om(\shfO,\shfO)\oplus {\mathcal
   H}om(\shfO(C),\shfO(C)).
\)

\begin{NB}
\textup{(4)} $\Ker d^0/\Ima d^{-1}$ is isomorphic to the cotangent bundle
to $\bp$ over $\Delta_{\bp}$.
\end{NB}
\end{Lemma}

\begin{Remark}
  After \propref{prop:diagonal} will be proved, it follows that the
  cohomology groups of the restriction of the complex $C^\bullet$ to
  the diagonal are $H^{-1} \cong K_{\bp}$, $H^0 \cong T^*\bp$, $H^1
  \cong \shfO$, since $C^\bullet\otimes\shfO_{\Delta}$ computes
  $\operatorname{Tor}^{\bp\times\bp}_\bullet(\shfO_\Delta,\shfO_\Delta)$.
  We prove this directly for $H^{\pm 1}$ in the above
  lemma. It is also possible to prove $H^0\cong T^*\bp$ directly.
\begin{NB}
  $\shfO_{\Delta_{\bp}}$ is quasi-isomorphic to $C^\bullet$. So
  $\shfO_\Delta\otimes \shfO_\Delta$ is quasi-isomorphic to
  $C^\bullet\otimes \shfO_\Delta$, which is nothing but the
  restriction of $C^\bullet$ to the diagonal $\Delta$.
  But from the definition of $\operatorname{Tor}$, this complex
  computes
  $\operatorname{Tor}^{\bp\times\bp}_\bullet(\shfO_\Delta,\shfO_\Delta)$.
  The cohomology groups are exterior powers of the cotangent bundles.
\end{NB}
\end{Remark}

\begin{proof}[Proof of \lemref{lem:diag}]
(1) From the definition $\varphi$ is the pull-back of the section
defined over $\proj^2\times\proj^2$. The latter section vanishes
on $\Delta_{\proj^2}$. Therefore $\varphi$ vanishes on
$(C\times C)\cup \Delta_{\bp}$.

\begin{NB}
The direct argument:
$\varphi$ vanishes if and only if $(s'z_0,z_1,z_2)$ is a multiple
of $(z_0',z',w')$. If $s'=0$ (and hence $z_0'\neq 0$), then the
multiplier must be $0$, hence we have $z_1 = z_2 = 0$, therefore
$([z_0,z_1,z_2],[z:w])\in C$. It means that the point is in
$C\times C$. Next suppose $s'\neq 0$. Then we have
\(
   [z_0:z_1:z_2] = [z_0'/s':z':w'] = [z_0':z_1':z_2'].
\)
As we have $s'\neq 0$, we must have $s\neq 0$ from this equation.
Therefore the second component $\proj^1$ in $\bp$ is determined by the
first factor, and hence the point is in the diagonal $\Delta_{\bp}$.
\end{NB}

(2) Suppose $d^{-1}$ is not injective. Take
\(
  a\oplus b\in (\shfO\boxtimes \shfO(C-2\linf))
   \oplus
   (\shfO(C-\linf)\boxtimes\shfO(-2\linf))
\)
from $\Ker d^{-1}$. Then $\varphi\lrcorner b = 0$. If $\varphi\neq 0$,
then we have $b=0$. It implies $z a = 0 = w a$, and hence $a = 0$.
Therefore $\varphi$ vanishes.
By (1) the point is in $(C\times C)\cup\Delta_{\bp}$. If it is in
$(C\times C)$, then we must have $z_0\neq 0$. Then looking at the
first and second row of $d^{-1}$, we must have $[z:w] = [z':w']$.
Conversely if the point is in the diagonal, it is clear that $d^{-1}$ is
not injective.
The second assertion is also obvious from the definition of $d^{-1}$.

(3) We have
\begin{equation*}
  (d^0)^\vee
  = 
   \begin{bmatrix}
    \id\boxtimes z_2 & s'\boxtimes w \\
    -\id\boxtimes z_1 & - s'\boxtimes z \\
    \chi & \psi
   \end{bmatrix}
   \colon
   \begin{matrix}
   \shfO\boxtimes\shfO \\ \oplus \\
   \shfO(-C)\boxtimes\shfO(C)
   \end{matrix} 
   \longrightarrow
\begin{matrix}
   \shfO^{\oplus 2}\boxtimes \shfO(\linf) 
   \\ \oplus \\ 
   \Qcal(-C)\boxtimes \shfO(\linf)
\end{matrix}.
\end{equation*}
Suppose that $0\neq \lambda\oplus \mu\in\shfO\boxtimes\shfO \oplus
\shfO(-C)\boxtimes \shfO(C)$ is in the kernel of $(d^0)^\vee$. We have
$\lambda \chi + \mu\psi = 0$, hence 
\begin{equation}\label{eq:1}
   (-\lambda z_0, \mu z, \mu w) \in \C(z_0',z',w'). 
\end{equation}
On the other hand, we have 
\begin{equation}
\label{eq:2}
  \lambda s + \mu s' = 0
\end{equation}
as $z_1 = sz$, $z_2 = sw$. 

Substituting \eqref{eq:2} to \eqref{eq:1}, we get
\(
   (\lambda z_0, \lambda z_1, \lambda z_2)
   \in \C (z_0', z_1', z_2').
\)
Therefore if $[z_0:z_1:z_2]\neq [z_0':z_1':z_2']$, we must have
$\lambda = 0$. Hence we have $\mu s' = 0$ from \eqref{eq:2}. On the
other hand, from \eqref{eq:1} and $\lambda = 0$, we have either
$\mu = 0$ or $z_0' = 0$. As $\{ s' = 0\}\cap \{ z_0' = 0\} =
\emptyset$, we have $\mu = 0$ anyway.

If $[z:w]\neq [z':w']$, we must have $\mu = 0$ and $\lambda z_0 = 0$
from \eqref{eq:1}. We also have $\lambda s = 0$ from \eqref{eq:2} and
$\mu = 0$. Therefore we must also have $\lambda = 0$ as $\{ s = 0\}\cap
\{ z_0 = 0\} = \emptyset$.

Finally suppose $[z_0:z_1:z_2] = [z_0':z_1':z_2']$ and $[z:w] =
[z':w']$. Then we have $\lambda + \mu = 0$ from either $z_0\neq 0$ or
$s\neq 0$.
\begin{NB}
(4) On the diagonal $\Delta_{\bp}$ we have
\begin{equation*}
   d^{-1} = \begin{bmatrix}
     z & z z_0 \\
     w & w z_0 \\
     0 & 0
     \end{bmatrix}, \quad
   d^0 = \begin{bmatrix}
    z_2 & - z_1 & \chi^{\vee}\\
    z_2 & - z_1 & \chi^{\vee}
   \end{bmatrix},
\end{equation*}
as $\psi = \chi$ on the diagonal. Outside $C$ (so $(z_1,z_2)\neq
(0,0)$), we can identify $\Ker d^0/\Ima d^{-1}$ with
$\Qcal^\vee(C-\linf)|_{\bp\setminus C}\cong \Qcal^\vee(-\linf)$ by the
natural projection to the factor $\Qcal^\vee(C-\linf)$.
Next consider a point
$(([1:z_1:z_2],[z:w]),([1:z_1:z_2],[z:w]))\in\Delta_{\bp}\setminus
\linf\times\linf$. We identify $\mathcal Q$ with $\shfO^{\oplus 2}$ as
above. Then we have
\begin{equation*}
   d^{-1} = 
   \begin{bmatrix}
     z & z
     \\
     w & w
     \\
     0 & 0
     \\
     0 & 0
   \end{bmatrix},
\qquad
   d^0 =
   \begin{bmatrix}
     z_2 & - z_1 & w & -z
     \\
     z_2 & - z_1 & w & -z
   \end{bmatrix}
\end{equation*}
Let us see that $\Ker d^0/\Ima d^{-1}$ is isomorphic to the cotangent
bundle from this description.
From the above, we identify $C^0/\Ima d^{-1}$ with
$\shfO(-C)\oplus \shfO(C)^{\oplus 2}$. And the induced
map $d^0$ on $C^0/\Ima d^{-1}$ is given by
\(
   \left[
   \begin{smallmatrix}
     s & w & -z
     \\
     s & w & -z
   \end{smallmatrix}\right].
\)
We describe $\widehat{\C}^2$ as the quotient of $\{ (z,w,s) \in \C^3
\mid (z,w)\neq (0,0) \}$ divided by $\C^*$ acting by $(z,w,s)\mapsto
(\lambda z,\lambda w,\lambda^{-1} s)$.
(In the description in \thmref{thm:blowupplane} we should identify
$(z,w,s) = (B_1,B_2,d)$.)
We then have a natural line bundle on $\widehat{\C}^2$ corresponding
to the $\C^*$-action, which can be identified with $\shfO(-C)$.
Moreover, $z$, $w$ can be considered its sections, and $s$ is a
section of its dual.
(Note that we are removing $\linf$.)
The tangent bundle of $\widehat{\C}^2$
is the quotient of
\(
  \shfO(-C)^{\oplus 2} \oplus \shfO(C)
\)
by the image 
\(
   \left[
     \begin{smallmatrix}
       z\\ w \\ - s
     \end{smallmatrix}
   \right]
\).

On the other hand, from the above description, we have
$\Ker (d^{-1})^\vee \cong \shfO(-C)^{\oplus 2} \oplus \shfO(C)$, where
\(
  \shfO(C) = \{ \left[
    \begin{smallmatrix}
      f \\ g
    \end{smallmatrix}
   \right] \in \shfO^{\oplus 2} \mid z f + w g = 0 \}.
\)
We can also identify $(d^0)^\vee$ with
\(
   \left[
     \begin{smallmatrix}
       z\\ w \\ - s
     \end{smallmatrix}
   \right]
\).

Now it is also clear that two isomorphisms to the cotangent bundle
glue on the intersection $\bp\setminus (\linf\cup C)$.
\end{NB}
\end{proof}

\begin{NB}
  Let us consider $\{ w'\neq 0\} \subset\bp$. Then we can identify
  $\mathcal Q$ with $\shfO^{\oplus 2}$ by
\begin{equation*}
  \begin{bmatrix}
      a \\ b \\ c
  \end{bmatrix}
  \bmod \Ima\iota
  \longmapsto
  \begin{bmatrix}
    a - b z_0' \\ c + bz'
  \end{bmatrix}.
\end{equation*}
The sections $\varphi$, $\psi$, $\chi$ are re-written as
\begin{gather*}
   \varphi = \left[
     \begin{smallmatrix}
      s'\boxtimes z_0 - z_0'\boxtimes z_2
      \\
      - \id\boxtimes z_1 + z'\boxtimes z_2
     \end{smallmatrix}
     \right],
\\
   \psi = \left[
       \begin{smallmatrix}
         -z_0'\boxtimes w \\ -\id\boxtimes z + z'\boxtimes w 
       \end{smallmatrix}
     \right],
\\
  \chi = \left[
       \begin{smallmatrix}
         -\operatorname{id}\boxtimes z_0 \\ 0
       \end{smallmatrix}
     \right].
\end{gather*}
\end{NB}

\begin{proof}[Proof of \propref{prop:diagonal}]
Let us first compute the alternating sum of Chern characters of
terms restricted to the slice $\{ x\}\times\bp$. The ranks and the
first Chern classes cancel. The alternating sum of the second Chern
characters is $1$.
\begin{NB}
\begin{equation*}
   \frac12 (C - 2\linf)^2 + \frac12 (-2\linf)^2 - \linf^2 - \linf^2 +
   \frac12 C^2
   = 1,
\end{equation*}
as $C^2 = -1$, $\linf^2 = 1$, $C\linf = 0$.

Let us restrict to the another slice $\bp\times\{x\}$. Then
$c_1(\Qcal) = \linf$, so the first Chern classes cancel. The second
Chern characters are
\begin{equation*}
   \frac12 (C-\linf)^2 - \ch_2(\Qcal) + c_1(\Qcal) \cdot C
   - C^2 + \frac12 C^2 = 1
\end{equation*}
as $\ch_2(\Qcal) = - 1/2$.
\end{NB}%
In particular (from the cancellation of ranks), the complex is exact
outside $\Delta$ by \lemref{lem:diag}.

Note the following: Suppose $A$, $B$ are flat over $T$, $f\colon A\to
B$ satisfies $f_x\colon A_x\to B_x$ is injective. Then $f$ is
injective and $B/A$ is flat over $T$.

From the above with \lemref{lem:diag}(2), we have $D := C^0/C^{-1}$ is
a family of torsion free sheaves on $\bp$ parametrized by $\bp$, and
$D_x$ is locally free over $\bp\setminus\{ x\}$. Next by
\lemref{lem:diag}(3) $D_x \to C^1_x$ is injective. So $D\to C^1$ is
injective and $C^1/D$ is a flat family of sheaves parametrized by
$\bp$.
Moreover $D_x\to C^1_x$ is an isomorphism on $\bp\setminus\{ x\}$. So
$C^1/D$ is a family of Artinian sheaves.
Finally we consider the induced homomorphism $C^1/D \to \shfO_{\Delta}$.
By \lemref{lem:diag}(3) $(C^1/D)_x\to\C_x$ is nonzero.
On the other hand, from the computation of the Chern classes, we see
$(C^1/D)_x\cong\C_x$. Thus $C^1/D\to \shfO_\Delta$ is an isomorphism.
\begin{NB}
The alternative argument:

The following proof is modeled on the corresponding statement for ALE
spaces in \cite[\S3]{KN}.

The only thing left doubt is whether the map
\begin{equation*}
   C^0 \to \Ker t
\end{equation*}
is surjective.

Let $\shfO_{\Delta,\bp\times\bp}$ be the localization of
$\shfO_{\bp\times\bp}$ along $\Delta = \Delta_\bp$, and $\mathfrak m$
be its maximal ideal. By Nakayama's lemma, it is enough to show that
\begin{equation*}
  C^0 \to (\Ker t)/\mathfrak m (\Ker t)
\end{equation*}
is surjective.

We interpret the right hand side in concrete terms.
Let $C^1|_\Delta$ be the restriction of the vector bundle $C^1$ to the
diagonal, considered as a locally free sheaf on $\Delta$. Let
$t|_\Delta\colon C^1|_\Delta\to\shfO_\Delta$ be the map given by
taking the difference of the trace.
Then $t$ factors as $C^1\to C^1|_\Delta\xrightarrow{t|_\Delta}
\shfO_\Delta$. Let $K$ be the kernel of $t|_\Delta$. Then we have the
surjective homomorphism
\begin{equation*}
   (\Ker t)/\mathfrak m (\Ker t) \to K
\end{equation*}
whose kernel is isomorphic to the the conormal bundle of $\Delta$ in
$\bp\times\bp$, which is identified with the cotangent bundle to
$\bp$. The isomorphism is given by taking the covariant derivative of
$f_1\in \Ker(\Ker t\to K)$ to the normal direction with respect to a
connection on $C^1$ and composing $t|_\Delta$. The result factors
through $\Ker[(\Ker t)/\mathfrak m(\Ker t)\to K]$ and is independent
of the choice of the connection.

From \lemref{lem:diag}(3), the composition of
\begin{equation*}
   C^0 \to (\Ker t)/\mathfrak m(\Ker t) \to K
\end{equation*}
is surjective. Therefore it is enough to show that
$\Ker[C^0\to K]\to T^*\bp$ is surjective.

The map factors as
\[ 
\xymatrix{
  \Ker[C^0\to K] \ar[rr] \ar@{->>}[rd] && T^*\bp, \\
   & 
   \frac{\Ker\left[C^0|_\Delta \xrightarrow{d^0|_\Delta} C^1|_\Delta\right]}
   {\Ima\left[C^{-1}|_\Delta \xrightarrow{d^{-1}|_\Delta} C^0|_\Delta\right]} 
   \ar[ru]^\varXi & \\
}
\] 
where $\varXi$ is given by $t|_\Delta \circ \nabla d^0$ for a
connection $\nabla$ on $(C^0)^*\otimes C^1$. It is enough to show that
$\varXi$ is surjective.

This is an identity map via \lemref{lem:diag}(4).%
\begin{NB2}
  Let us consider the map over $\Delta\setminus\linf\times\linf$:
  Let us compute the derivative of $d^0$ along
  $\Delta$. Over $w\neq 0$ we take the coordinates $(z,z_2)$ (and put
  $w=1$). Then
\begin{equation*}
    d^0
   =
   \begin{bmatrix}
     z_2 & - z z_2 & 1 & -z'
     \\
     z_2' & - z_2' z & 1 & - z
   \end{bmatrix}
\end{equation*}
and
\begin{equation*}
   t|_\Delta \circ \nabla d^0
   =
   \begin{bmatrix}
     dz_2 - dz'_2 & -z(dz_2 - dz'_2) & 0 & dz-dz'
   \end{bmatrix}.
\end{equation*}
Next consider $\Delta\cap(\linf\times\linf)$.%
We take coordinates $(z_0,z)$ around $([0:0:1],[0:1])$ given by
$([z_0:z:1],[z:1])$. Then
\begin{equation*}
   d^0 =
   \begin{bmatrix}
     1 & -z & - z_0 & 0 \\
     1 & -z & -z'_0 & -z + z'
   \end{bmatrix}.
\end{equation*}
Therefore
\begin{equation*}
   t|_\Delta\circ \nabla d^0 =
   \begin{bmatrix}
     0 & 0 & d z'_0 - d z_0 & dz' - dz \\
   \end{bmatrix}.
\end{equation*}
\end{NB2}
\end{NB}
\end{proof}

\section{Analysis of stability}\label{sec:analysis}

We fix $\zeta=(\zeta_1,\zeta_2)$ as in \eqref{eq:ass}. 

\subsection{Necessary conditions for stability}\label{subsec:nec}

We shall use the language of the representation theory of a quiver
with relation. The data $(B_1,B_2,d,i,j)$ with the underlying vector
spaces $V_0$, $V_1$, $W$, as in \thmref{thm:King} is a representation
of a quiver with relation. But as $W$ plays a different role from
$V_0$, $V_1$, it is more natural to construct a new quiver following
\cite{CB}.

\begin{figure}
  \centering
\begin{equation*}
\xymatrix{
0 \ar[rrrr] &&&& \ar@<-1ex>[llll] \ar@<-2ex>[llll] \ar@/^1pc/[lldd]_\ddots
\ar@<-1ex>[llll] \ar@/_1pc/[lldd]  1
\\
&&\\
&& \ar@/^1pc/[lluu]_\iddots
  \ar@/_1pc/[lluu]  \infty & \\
}
\end{equation*}
  \caption{new quiver}
  \label{fig:newquiver}
\end{figure}
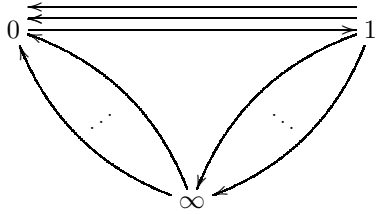

We fix the vector space $W$ with $\dim W = r$.
We define a new quiver with three vertexes $0$, $1$, $\infty$. We
write two arrows from $1$ to $0$ corresponding to the data $B_1$,
$B_2$, and one arrow from $0$ to $1$ corresponding to the data $d$.
Instead of writing one arrow from $\infty$ to $0$, we write
$r$-arrows. Similarly we write $r$-arrows from $1$ to $\infty$.
And instead of putting $W$ at $\infty$, we replace it the
one dimensional space $\C$ on $\infty$. We denote it by $V_\infty$.
It means that instead of considering the homomorphism
$i$ from $W$ to $V_0$, we take $r$-homomorphisms
$i_1$,$i_2$,\dots, $i_{r}$ from $V_\infty$ to $V_0$ by taking a base of $W$.
(See Figure~\ref{fig:newquiver}.)

We consider the full subcategory of the abelian category of
representations of the new quiver with the relation, consisting of
representations such that $\dim V_\infty = 0$ or $1$. An object can be
considered as a representation of the original quiver with $\dim W =
0$, or $\dim W = r$, according to $\dim V_\infty = 0$ or $1$. Note
that we do not allow a representation of the original quiver with
$\dim W\neq r,0$.
For objects $X$ and $Y$, let $\Hom(X,Y)$ denote the space of morphisms
from $X$ to $Y$. 

\begin{NB}
Let $X = (B_1,B_2,d,i,j)$ be the data satisfying 
\(
   \mu(B_1, B_2,d,i,j) = 0.
\)
Suppose that we have another data $Y = (B_1',B_2',d')$ satisfying
\(
   \mu(B_1',B_2', d') = B'_1 d' B'_2 - B'_2 d' B'_1 = 0
\)
defined for another pair of vector spaces $V_0'$, $V_1'$ with $W' =
0$. A {\it homomorphism\/} from $X = (B_1,B_2,d,i,j)$ to $Y = (B_1',B_2',d')$
is a pair of homomorphisms $\xi_0\in \Hom(V_0,V'_0)$,
$\xi_1\in\Hom(V_1,V_1')$ such that all diagrams
\begin{equation*}
  \begin{CD}
     V_1 @>B_{\alpha}>> V_0
     \\
     @V{\xi_1}VV @VV{\xi_0}V
     \\
     V'_1 @>B'_{\alpha}>> V'_0
  \end{CD}
  \quad (\alpha=1,2), \qquad
  \begin{CD}
     V_0 @>d>> V_1
     \\
     @V{\xi_0}VV @VV{\xi_1}V
     \\
     V'_0 @>d'>> V'_1
  \end{CD}
\end{equation*}
are commutative. We denote the space of homomorphisms by
$\Hom(X,Y)$. It is a complex vector space. Similarly we define a
homomorphism $\eta$ from $Y$ to $X$. 
\end{NB}

Recall $C_m = (B_1,B_2,0)$ denote the data given in \thmref{thm:W=0}.

\begin{Proposition}\label{prop:ss}
Suppose $X = (B_1,B_2,d,i,j)$ satisfies
\(
   \mu(B_1,B_2,\allowbreak d,i,j) = 0.
\)
If $X$ is $\zeta$-semistable, then the following holds:
\begin{enumerate}
\item If $m\zeta_0 + (m+1)\zeta_1 > 0$, then
$\Hom(C_m, X) = 0$.
\item If $m\zeta_0 + (m+1)\zeta_1 < 0$, then
$\Hom(X,C_m) = 0$.
\end{enumerate}
If $X$ is $\zeta$-stable, then we also have
\begin{enumerate}
\setcounter{enumi}{2}
\item If $m\zeta_0 + (m+1)\zeta_1 = 0$, then
$\Hom(X,C_m) = 0 = \Hom(C_m,X)$.
\end{enumerate}
\end{Proposition}

\begin{proof}
(1) Suppose $X = (B_1,B_2,d,i,j)$ is $\zeta$-semistable, and take
$m$ with $m\zeta_0 + (m+1)\zeta_1 > 0$ and a homomorphism
$\xi\colon C_m\to X$. Consider $S_0 = \Ima\xi_0$, $S_1 =
\Ima\xi_1$. They satisfy the assumption in the stability condition 
for $X$ in \defref{def:stable}. By the $\zeta$-semistability of $X$ we
have
\begin{equation*}
   \zeta_0 \dim S_0 + \zeta_1 \dim S_1 \le 0.
\end{equation*}
On the other hand, $S'_0 = \Ker\xi_0$, $S'_1 = \Ker\xi_1$ satisfy the
assumption in the stability condition for $C_m$ in
\defref{def:stable}. As $C_m$ is $(-1, m/(m+1))$-stable by
\lemref{lem:C_mstable} we have
\begin{equation*}
   - \dim S'_0 + \frac{m}{m+1} \dim S'_1 \le 0,
\end{equation*}
or equivalently
\begin{equation*}
   - \dim S_0 + \frac{m}{m+1} \dim S_1 \ge 0.
\end{equation*}
But the two inequalities contradict with the assumption $m\zeta_0 +
(m+1)\zeta_1 > 0$ unless $\dim S_0 = \dim S_1 = 0$.%
\begin{NB}
We have $-\frac{\zeta_1}{\zeta_0} > \frac{m}{m+1}$, hence
\begin{equation*}
  0 \le -\dim S_0 + \frac{m}{m+1} \dim S_1 \le
  -\dim S_0 - \frac{\zeta_1}{\zeta_0}\dim S_1 \le 0.
\end{equation*}
The second inequality is strict unless $\dim S_1 = 0$. Therefore we
have $\dim S_1 = 0$. Substituting it back to two inequalities, we get
$\dim S_0 = 0$.
\end{NB}

The statement (2) can be proved in the same way.

(3) The above proof works until the last stage under $m\zeta_0 +
(m+1)\zeta_1 = 0$. The first inequality and the last inequality are in
the opposite direction, so we have equalities in all the three
inequalities. As $X$ is further $\zeta$-stable, we have $S_0 = 0$,
$S_1 = 0$.
\end{proof}

\subsection{Sufficient conditions for stability}

The purpose of this section is to show a partial converse to the
previous proposition:

\begin{Proposition}\label{prop:ss2}
  Suppose that $X = (B_1,B_2,d,i,j)$ satisfies $\mu(B_1,\allowbreak
  B_2,d,i,j) = 0$ and the condition \textup{(S2)} in
  \thmref{thm:King}. If $X$ satisfies \textup{(1)} and \textup{(2)} in
  \propref{prop:ss} for a given $\zeta$, then $X$ is
  $\zeta$-semistable.
\end{Proposition}

Note that (S2) is a necessary condition for $\zeta$-semistability by
\lemref{lem:S2}. Hence $X$ is $\zeta$-semistable if and
only if $X$ satisfies (S2), \textup{(1)} and \text{(2)} in
\propref{prop:ss}. A meaning of the condition (S2) will be clarified
in \lemref{lem:King}.

For the proof we use the Harder-Narasimhan and Jordan-H\"older
filtration of a representation of a quiver.
For this purpose, we slightly need to modify the stability condition,
which is suitable for representations for the new quiver. We fix data
$(B_1,B_2,d,i,j)$ defined for $V_0$, $V_1$, $W$, and consider it as a
representation $X$ of the new quiver as before. We define
$\zeta_\infty$ by 
\( 
   \zeta_\infty = - \zeta_0 \dim V_0 - \zeta_1 \dim V_1.
\) 
For another nonzero representation $Y$ of the new quiver with the
underlying vector space $V_0'\oplus V_1'\oplus V_\infty'$, we define
the {\it slope\/} by
\begin{equation*}
    \theta(Y) \defeq \frac{\zeta_0 \dim V_0' + \zeta_1 \dim V_1' +
      \zeta_\infty \dim V_\infty'}{\dim V_0'+\dim V_1'+\dim V_\infty'}.
\end{equation*}
We only consider the case $\dim V_\infty' = 0$ or $1$ as before. We
say $Y$ is {\it $\theta$-semistable\/} if we have
\begin{equation*}
    \theta(S) \le \theta(Y)
\end{equation*}
for any subrepresentation $0\neq S$ of $Y$. We say $Y$ is {\it
  $\theta$-stable\/} if the inequality is strict unless $S = Y$. If
$(\dim V_0',\dim V'_1,\dim V'_\infty) = \linebreak[4]
(\dim V_0,\dim V_1, 1)$ (for
example, when $Y = X$), then $\theta(Y) = 0$. In this case,
$\theta$-(semi)stability is equivalent to $\zeta$-(semi)stability.
In fact, if a subrepresentation $S$ has $S_\infty = 0$, then
$\theta(S) \le 0$ is equivalent to $\zeta_0\dim S_0 + \zeta_1\dim
S_1\le 0$. If a subrepresentation $T$ has $T_\infty = \C$, then
$\theta(T)\le 0$ is equivalent to $\zeta_0\codim T_0 + \zeta_1\codim
T_1\ge 0$.

When $(\dim V_0',\dim V'_1,\dim V'_\infty)$ is not equal to $(\dim
V_0,\dim V_1, 1)$, we define a new parameter $\zeta' =
(\zeta'_0,\zeta'_1,\zeta'_\infty)$ by
\begin{equation}\label{eq:newdata}
    (\zeta'_0,\zeta'_1,\zeta'_\infty)
    = (\zeta_0,\zeta_1,\zeta_\infty) 
        - \theta(Y)(1,1,1).
\end{equation}%
\begin{NB}
In particular, we have
\(
    \zeta'_0 \dim V'_0 + \zeta'_1 \dim V'_1 + \zeta'_\infty \dim
    V'_\infty = 0.
\)  
\end{NB}
Then the $\theta$-(semi)stability of $Y$ is equivalent to
$(\zeta'_0,\zeta'_1)$-(semi)stability of $Y$, as
\[
    \frac{\zeta'_0 \dim S_0 + \zeta'_1 \dim S_1 
    + \zeta'_\infty \dim S_\infty}{\dim S_0 + \dim S_1 + \dim S_\infty}
    = \theta(S) - \theta(Y).
\]
Suppose $\dim V'_\infty = 0$. Then
\begin{equation}\label{eq:ineq}
\begin{split}
  & \zeta_0' + \zeta_1' = \frac1{\dim V'_0 + \dim V'_1}
  (\zeta_0 - \zeta_1)(\dim V'_1 - \dim V'_0),
\\
  & \zeta'_0 = \frac{\dim V'_1}{\dim V'_0+\dim V'_\infty}
  (\zeta_0 - \zeta_1).
\end{split}
\end{equation}

\begin{Theorem}\cite{Rudakov}\label{thm:Rudakov}
\textup{(1)} A representation $Y$ has a Harder-Narasimhan filtration:
\begin{equation*}
   Y = Y^0 \supset Y^1 \supset \cdots \supset Y^N \supset Y^{N+1} = 0
\end{equation*}
such that $Y^k/Y^{k+1}$ is $\theta$-semistable for $k=0,1,\dots,N$ and
$\theta(Y^0/Y^1) < \theta(Y^1/Y^2) < \cdots < \theta(Y^N/Y^{N+1})$.

\textup{(2)} A $\theta$-semistable representation $Y$ has a
Jordan-H\"older filtration:
\begin{equation*}
   Y = Y^0 \supset Y^1 \supset \cdots \supset Y^N \supset Y^{N+1} = 0
\end{equation*}
such that $Y^k/Y^{k+1}$ is $\theta$-stable for $k=0,1,\dots,N$ and
$\theta(Y^0/Y^1) = \theta(Y^1/Y^2) = \cdots = \theta(Y^N/Y^{N+1})$.
\end{Theorem}

\begin{proof}[Proof of \propref{prop:ss2}]
Suppose that $X$ is not $\zeta$-semistable. Since
$\zeta$-semistability is equivalent to the $\theta$-semistability in
this case, $X$ is not $\theta$-semistable. So we take a
Harder-Narasimhan filtration $X = X^0 \supset \cdots \supset X^N
\supset X^{N+1} = 0$ as in \thmref{thm:Rudakov}. We have $N > 0$ by
the assumption. If
$\theta(X^0/X^1) \ge 0$, then $\theta(X^k/X^{k+1}) > 0$ for
$k=1,\dots,N$. It contradicts with $\theta(X) = 0$. Therefore
$\theta(X^0/X^1) < 0$. Similarly we must have $\theta(X^N/X^{N+1}) > 0$.

Suppose $(X^0/X^1)_\infty = 0$. Then $T_0 = (X^1)_0$, $T_1 = (X^1)_1$
satisfy the condition in (S2) in \thmref{thm:King}. As we assume (S2),
we have $\codim T_0 < \codim T_1$, i.e., $\dim (X^0/X^1)_0 < \dim
(X^0/X^1)_1$. If we define the new parameter $\zeta' =
(\zeta'_0,\zeta'_1,\zeta'_\infty)$ by the formula \eqref{eq:newdata}
for the representation $X^0/X^1$, we have
\(
   \zeta'_0 + \zeta'_1 < 0,
\)
and
\(
   \zeta'_0 < 0
\)
by \eqref{eq:ineq} and the assumption \eqref{eq:ass}. Therefore
\eqref{eq:ass} is satisfied for $(\zeta'_0,\zeta'_1)$.
If $X^0/X^1$ is $\theta$-stable, or equivalently
$(\zeta'_0,\zeta'_1)$-stable, then by \thmref{thm:W=0} we have
$X^0/X^1 = C_m$ for some $m$. But $\theta(X^0/X^1) = \theta(C_m) < 0$
means $m \zeta_0 + (m+1)\zeta_1 < 0$, hence $X = X^0$ cannot have
$C_m$ as a quotient by our assumption (2) in \propref{prop:ss}. This
is a contradiction.
Even if $X^0/X^1$ is not $\theta$-stable, we take a Jordan-H\"older
filtration of $X^0/X^1$, in particular $X$ has $C_m$ as a quotient
with $\theta(C_m) < 0$. This is again a contradiction.
Therefore we must have $(X^0/X^1)_\infty = \C$.

Next assume $(X^N/X^{N+1})_\infty = X^N_\infty = 0$. Then we have
\begin{equation*}
  0 < \theta(X^N) = \frac{\zeta_0 \dim X^N_0 + \zeta_1 \dim X^N_1}
  {\dim X^N_0 + \dim X^N_1}.
\end{equation*}
By the assumption \eqref{eq:ass} we must have $\dim X^N_0 < \dim
X^N_1$. The rest of the argument goes along the same line as above. We
must have $C_m \subset X$ with $m\zeta_0 + (m+1)\zeta_1 > 0$, but it
contradicts with our assumption (1) in \propref{prop:ss}. Therefore we
must have $(X^N/X^{N+1})_\infty = \C$, but as $X_\infty$ is
$1$-dimensional, it cannot be possible to have both
$(X^0/X^1)_\infty = \C$ and $(X^N/X^{N+1})_\infty = \C$ unless
$N=0$. It means that $X$ is $\zeta$-semistable.
\end{proof}

\subsection{Walls and chambers}\label{subsec:walldef}
Take $\zeta = (\zeta_0,\zeta_1)\in\R^2$ as in \eqref{eq:ass}. Suppose
$\zeta$-semistability is not equivalent to the $\zeta$-stability. Take
data $X = (B_1,B_2,d,i,j)\in\mu^{-1}(0)$ which is $\zeta$-semistable,
but not $\zeta$-stable. We consider $X$ as a representation of the new
quiver as above. Then by \thmref{thm:Rudakov} we have a
Jordan-H\"older filtration $X = X^0 \supset \cdots \supset X^N \supset
X^{N+1} = 0$ as in \thmref{thm:Rudakov}. We have $N > 0$ by the
assumption. As $\dim X_\infty = 1$, we have either $\dim
(X^0/X^1)_\infty = 0$ or $\dim (X^N/X^{N+1})_\infty = 0$. In
particular, we have a $\zeta$-stable representation $Y$ with $Y_\infty
= 0$. By \thmref{thm:W=0} then $Y = C_m$ for some $m\in\Z_{\ge 0}$ and
we have $m\zeta_0 + (m+1)\zeta_1 = 0$. Thus this is a possible wall in
the space of parameters $\zeta$ as in \eqref{eq:ass}.

A {\it wall of type $(r, k, n)$\/} is, by definition, a hyperplane $\{
m\zeta_0+(m+1)\zeta_1 = 0\}$ such that there is $X\in\mu^{-1}(0)$ as
above with the data $r$, $k$, $n$. For fixed $r,k,n$, walls of type
$(r,k,n)$ is finite. A connected component of the complement of walls
of type $(r,k,n)$ is called a {\it chamber of type $(r,k,n)$}. If
there is no fear of confusion, we simply say a chamber.
If $\zeta$ stays in a chamber, the corresponding stability condition
remains equivalent. When we cross a wall, then the stability condition
and the moduli space are changed. 

\section{Construction : from ADHM data to a sheaf}\label{sec:const}

In this section we construct a coherent sheaf $E$ on $\bp$ together
with a framing $\Phi$ such that $(E(-mC),\Phi)$ is stable perverse coherent
from a $\zeta$-stable data $X = (B_1,B_2,d,i,j)$.

\subsection{Surjectivity of $\beta$}

\begin{Lemma}[\protect{\cite[4.1.3]{King}}]\label{lem:King}
$X = (B_1,B_2,d,i,j)$ satisfies the condition \textup{(S2)} in
\thmref{thm:King} if and only if $\beta$ in the complex \eqref{eq:cpx}
is surjective at every point in $\bp$.
\end{Lemma}

\begin{proof}
For the sake of completeness, we give a proof. The proof is
essentially the same as the original one, but our use of
\thmref{thm:blowupplane} somewhat clarifies the role of
\cite[Lemma~4.1.1, Remark~4.1.2]{King}.

Let us first show the `only if' part:
Assume $\beta$ is not surjective at a point $([z_0:z_1:z_2],[z:w])$.
As $\beta$ is surjective along $z_0 = 0$, we may assume $z_0 = 1$. Then
there exists $v_0\in V_0^*$, $v_1\in V_1^*$ such that
\begin{equation}\label{eq:not_surj}
\begin{gathered}
  z_2 v_0 + w d^* v_1 = 0, \quad
  z_1 v_0 + z d^* v_1 = 0, \quad
  B_2^* v_0 + w v_1 = 0, \\
  B_1^* v_0 + z v_1 = 0, \quad
  i^* v_0 = 0
\end{gathered}
\end{equation}
and $(v_0,v_1)\neq (0,0)$. Then $(S_0, S_1) \defeq (\C v_0, \C v_1)$
satisfies $B_\alpha^* (S_0)\subset S_1$, $d^*(S_1)\subset S_0$,
$i^*(S_0) = 0$. Moreover we have $v_0\neq 0$, as $v_0 = 0$ implies
$v_1 = 0$.  Therefore $\dim S_0 = 1$ and $\dim S_1 = 0$ or $1$. Hence
the annihilators $T_0 \defeq S_0^\perp\subset V_0$, $T_1\defeq
S_1^\perp \subset V_1$ violate the condition (S2).

We prove the `if' part: Assume $X = (B_1,B_2,d,i,j)$ violates the
condition (S2). We consider $(B_1^*,B_2^*,d^*,i^*,j^*)$, then we have
a pair of subspaces $S_0\subset V_0^*$, $S_1\subset V_1^*$ such that
$B_\alpha^*(S_0)\subset S_1$ ($\alpha=1,2$), $d^*(S_1)\subset S_0$,
$i^*(S_0) = 0$, and $\dim S_1\le \dim S_0$ with $(S_0,S_1)\neq (0,0)$.
Using the idea of the proof of the Harder-Narasimhan filtration, we
take $(S_0,S_1)$ among such pairs of subspaces, so that
\begin{equation*}
   \theta(S_0,S_1) \defeq \frac{\dim S_0 - \dim S_1}{\dim S_0 + \dim S_1}
\end{equation*}
takes the maximum. In particular, $\theta(S_0,S_1)\ge 0$. Then the
restriction of $(B_1^*,B_2^*,d^*)$ to $(S_0,S_1)$ is
$\theta$-semistable. We then consider a Jordan-H\"older filtration
with respect to $\theta$, we may further assume that the restriction
of $(B_1^*,B_2^*,d^*)$ to $(S_0,S_1)$ is $\theta$-stable. If we
define a new parameter $(\zeta'_0,\zeta'_1)$ by
\begin{equation*}
   (\zeta'_0,\zeta'_1) = (1,-1) 
   - \theta(S_0,S_1)(1,1),
\end{equation*}
we have $\zeta'_0 \dim S_0 + \zeta'_1 \dim S_1 = 0$ and
$\theta$-stability of $(B_1^*,B_2^*,d^*)|_{(S_0,S_1)}$ is equivalent
to its $(\zeta'_0,\zeta'_1)$-stability in \defref{def:stableW=0}. We
have $\zeta'_0 + \zeta'_1 \le 0$ and $\zeta'_1 < 0$. Remarking that we
change the numbering of the vertexes $0$, $1$ as we are considering
the transposes $B_1^*$, $B_2^*$, $d^*$, the parameter
$(\zeta'_0,\zeta'_1)$ satisfies either \eqref{eq:ass} or
\eqref{eq:ass'}. Therefore the restriction of $(B_1^*, B_2^*, d^*)$ to
$(S_0,S_1)$ is classified as in \thmref{thm:W=0} and
\thmref{thm:blowupplane}. In the latter case, there is a corresponding
point $(z_1,z_2,[z:w])\in \widehat{\mathbb C}^2$. The surjectivity of
$\beta$ is violated at the point. In the former case, there exists
$0\neq v_0\in S_0$ such that $(w B_1^* - z B_2^*)v_0 = 0$ as $\dim S_0
> \dim S_1$ for any pair $(z,w)$. Therefore we can find vectors
$v_0\in V_0^*$, $v_1\in V_1^*$ as in \eqref{eq:not_surj} for any point
on the exceptional curve $C = \{ (0,0,[z:w])\}$. Again the
surjectivity of $\beta$ is violated.
\end{proof}

\subsection{The definition of a framed sheaf}\label{subsec:defSheaf}

Suppose that $(B_1,B_2,d,i,j)$ satisfies the condition (S2). Then
$\beta$ in \eqref{eq:cpx} is surjective at every point by
\lemref{lem:King}. Let us consider the restriction of \eqref{eq:cpx}
to the line $\linf$ at infinity.
\begin{NB}
\begin{equation*}
      \left.\alpha\right|_{\linf} = 
   \begin{bmatrix}
     z & 0           \\
     w & 0           \\
     0 & z_1  \\
     0 & z_2  \\
     0 & 0
   \end{bmatrix},
\qquad
   \left.\beta\right|_{\linf} =
   \begin{bmatrix}
    z_2 & - z_1 & 0       & 0        & 0\\
    dw  & - d z & w       & -z       & 0
   \end{bmatrix}.
\end{equation*}
\end{NB}%
Then we can check that $\left.\alpha\right|_{\linf}$ is injective
and 
\(
   (\left.\Ker \beta\right|_{\linf}) / (\Ima \left.\alpha\right|_{\linf})
\)
is isomorphic to $W\otimes\shfO_{\linf}$ via the natural inclusion
of $W\otimes\shfO_{\linf}$ to the middle term of \eqref{eq:cpx}
restricted to $\linf$. Hence we have
\begin{enumerate}
\item $\alpha$ is injective as a sheaf homomorphism, and hence
  the complex \eqref{eq:cpx} is quasi-isomorphic to a coherent
  sheaf $E = \Ker\beta/\Ima \alpha$,
\item the sheaf $E$ has the canonical framing 
  $\varphi\colon E|_{\linf} \xrightarrow{\cong} W\otimes \shfO_{\linf}$.
\end{enumerate}

Furthermore the injectivity of $\alpha$ fails along a subvariety of
$\bp$, therefore $\alpha$ is injective except possibly along $C$ and
finite points in $\widehat{\mathbb C}^2\setminus C$. In particular,
$E$ is torsion free outside $C$ (see e.g.\ the argument in
\cite[\S5(i)]{Na:ADHM}).

\subsection{The sheaf corresponding to $C_m$}

For a later purpose we notice that this construction works in the case
$W = 0$, provided the assumption (S2) is satisfied. For example, the
data $C_m$ in \thmref{thm:W=0} gives a sheaf: $C_m$ satisfies
(S2) by \lemref{lem:S2} as it is $(-1,m/(m+1))$-stable. We can
explicitly describe it:

\begin{Proposition}\label{prop:Cm}
The sheaf corresponding to $C_m$ is 
\(
   \shfO_C(-m-1).
\)
\end{Proposition}

\begin{proof}
We use the resolution of the diagonal $C^\bullet$
in \eqref{eq:diagonal}. Our proof is a prototype of the argument used in
\secref{sec:inverse}. We consider 
\[
   Rp_{2*} \left(p_1^* \shfO_C(-m-1)\otimes C^\bullet\right),
\]
where $p_1$, $p_2$ are the projection to the first and second factors
of $\bp\times\bp$ respectively. Since $C^\bullet$ is quasi-isomorphic
to the diagonal, this is isomorphic to $\shfO_C(-m-1)$.

Next let us take the tensor product of $p_1^*\shfO_C(-m-1)$ with
$C^\bullet$ in \eqref{eq:diagonal} and take derived push-forwards.
Note the cohomology groups $H^i(\shfO_C(-m-1))$ and
$H^i(\shfO_C(-m-1)\otimes\shfO(C)) = H^i(\shfO_C(-m-2))$ vanishes for
$i = 0$, and have dimensions $m$ and $m+1$ respectively for $i = 1$
thanks to the assumption $m \ge 0$.
Let $V_0 = H^1(\shfO_C(-m-1))$, $V_1 =
H^1(\shfO_C(-m-1)\otimes\shfO(C))$.
We consider the homomorphisms on cohomology groups induced from $d_1$,
$d_2$ in \eqref{eq:diagonal}. As $z'_1 = z'_2 = s = 0$ on $C =
\operatorname{Supp}\shfO_C(-m-1)$, the only nontrivial homomorphisms
are those $V_1\to V_0$ induced from $z'$, $w'$. Under a suitable
choice of bases of $V_1$, $V_0$, the homomorphisms are given by
$B_1$, $B_2$ in (\ref{eq:Kronecker}b).%
\begin{NB}
  The linear maps $B_1$, $B_2$ is multiplication by $z$, $w$, and
  via the Serre duality, they are the transpose of the map
\begin{equation*}
   z, w \colon 
   V_0^* = H^0(\shfO_C(m-1)) = 
   \operatorname{Span}(z^{m-1}, z^{m-2}w, \dots, w^{m-1})
   \to
   V_1^* = H^0(\shfO_C(m)) = 
   \operatorname{Span}(z^{m}, z^{m-1}w, \dots, w^{m}).
\end{equation*}
Now the claim is obvious.
\end{NB}
Then we see that the homomorphisms induced from $d_2$, $d_1$ are
exactly equal to $\alpha$, $\beta$ in \eqref{eq:cpx} for the data
$C_m$. Therefore we conclude that $\Ker\beta/\Ima\alpha$ for $C_m$ is
isomorphic to $\shfO_C(-m-1)$.
\end{proof}

\begin{Remark}
  The framing of $C_m$ is implicit as $E|_{\linf} \cong 0$.
\end{Remark}

\subsection{Homomorphisms}

Suppose that two framed coherent sheaves $(E,\Phi)$ and $(E',\Phi')$
are obtained by the construction in \subsecref{subsec:defSheaf}. Let
us denote the data for $E$ by $V_0$, $V_1$, $B_1$, etc, and by $V'_0$,
$V'_1$, $B'_1$ for $E'$. Let us denote the corresponding representations
of the quiver (with relations) by $X$ and $X'$ respectively.
Let $\Hom(X,X')$ be the space of homomorphisms from $X$
to $X'$, considered as representations of the quiver.
Here the quiver is as in \thmref{thm:King} and is not the
{\it new quiver\/}, used in \secref{sec:analysis}.
But the difference between the old and new quivers is not relevant, as
we are interested only in the cases when one of $W$ or $W'$ is $0$ (as
in \secref{sec:analysis}), or $W = W'$. In the latter case, we are only
interested in homomorphisms which are identity on $W = W'$. Anyway, we
use the same notation $\Hom$ as in \secref{sec:analysis}, as there is
no fear of confusion.

\begin{Proposition}\label{prop:hom}
  Let us consider the space of homomorphisms $\xi\colon E\to E'$. It
  is isomorphic to the space of homomorphisms $\Hom(X,X')$, i.e.,
  pairs of linear maps $\xi_0\colon V_0\to V'_0$, $\xi_1\colon V_1\to
  V'_1$, $b_{WW}\colon W\to W'$ such that
  \begin{equation}\label{eq:commute}
    \begin{gathered}
     \xi_0 B_\alpha = B'_\alpha \xi_1 \quad (\alpha=1,2), \qquad
     \xi_1 d = d' \xi_0, \\
     \xi_0 i = i'b_{WW}, \qquad b_{WW} j = j' \xi_1.
    \end{gathered}
  \end{equation}
\end{Proposition}

Here we understood $i = 0 = j$ (resp.\ $i' = 0 = j'$) if
$W = 0$ (resp.\ $W' = 0$).

\begin{NB}
We only consider the two cases
\begin{aenume}
\item $W = 0$ or $W' = 0$,
\item $\dim W = \dim W'$.
\end{aenume}
In the second case, we also fix an isomorphism $W \cong W'$ of vector
spaces, and identify $W$ with $W'$. We also have an induced
identification $E|_{\linf}\cong E'|_{\linf}$ via the framings.
\end{NB}

For the proof, we need a terminology from the theory of monads: A {\it
  monad\/} is a complex of vector bundles
\begin{equation*}
   0 \to \mathcal A \xrightarrow{\alpha} \mathcal B
   \xrightarrow{\beta}
   \mathcal C \to 0
\end{equation*}
such that $\alpha$ is injective and $\beta$ is surjective as sheaf
homomorphisms, so that it is quasi-isomorphic to a coherent sheaf
$E = \Ker\beta/\Ima\alpha$. A morphism from one monad
\( 
  M\colon \mathcal A \xrightarrow{\alpha} \mathcal B
   \xrightarrow{\beta} \mathcal C
\)
to another
\(
  M'\colon \mathcal A' \xrightarrow{\alpha'} \mathcal B'
   \xrightarrow{\beta'} \mathcal C'
\)
is a morphism of complexes, i.e., the diagram
\begin{equation}\label{eq:monadHom}
  \begin{CD}
   \mathcal A @>{\alpha}>> \mathcal B
   @>{\beta}>> \mathcal C
\\
   @VVaV @VVbV @VVcV
\\
    \mathcal A' @>>{\alpha'}> \mathcal B'
   @>>{\beta'}> \mathcal C'.
  \end{CD}
\end{equation}
It is determined by $b$, as long as it satisfies $b(\Ima\alpha)
\subset\Ima\alpha'$, $b(\Ker\beta)\subset\Ker\beta'$. Then we have an
induced map $H(b)\colon E\to E'$. Let $\Hom_{\mathrm{mon}}(M,M')$
denote the set of morphisms from $M$ to $M'$.

\begin{Lemma}[\protect{\cite[2.2.1]{King}}]\label{lem:monad}
  Let $M$, $M'$ be monads as above, and $E$, $E'$ be corresponding
  sheaves. Suppose that
\begin{equation*}
    \Hom(\mathcal B, \mathcal A'),\quad 
    \Hom(\mathcal C, \mathcal B'),\quad
    \Ext^1(\mathcal C, \mathcal B'),\quad
    \Ext^1(\mathcal B, \mathcal A'),\quad
    \Ext^2(\mathcal C, \mathcal A')
\end{equation*}
all vanish. Then the above assignment
\(
   H\colon \Hom_{\mathrm{mon}}(M,M') \to \Hom(E,E')
\)
is surjective, and the kernel is isomorphic to $\Ext^1(\mathcal C,
\mathcal A')$.
\end{Lemma}

\begin{proof}
We give a proof following \cite[Ch.~II, Lem.~4.1.3]{OSS}.

Let us suppose a homomorphism $\varphi\colon E\to E'$ is given. We
claim that there exists a homomorphism
$\psi\colon\Ker b\to \Ker b'$ completing the commutative diagram
(with exact rows)
\begin{equation*}
  \begin{CD}
    0 @>>> \mathcal A @>>> \Ker \beta @>>> E @>>> 0
\\
    @. @VVV @VV{\psi}V @VV{\varphi}V @.
\\
    0 @>>> \mathcal A' @>>> \Ker \beta' @>>> E' @>>> 0,
  \end{CD}
\end{equation*}
and it is unique up to $\Hom(\Ker \beta,\mathcal A')$.

First we prove
\begin{equation}\label{eq:coh1}
  \Hom(\Ker \beta,\mathcal A') \cong \Ext^1(\mathcal C,\mathcal A'),
  \qquad
  \Ext^1(\Ker \beta,\mathcal A') = 0
\end{equation}
We apply the derived functor ${\mathbf R}\Hom(\bullet,\mathcal A')$ to
the exact sequence
\begin{equation}\label{eq:exact1}
   0 \to \Ker \beta\to \mathcal B \to \mathcal C \to 0.
\end{equation}
Then \eqref{eq:coh1} follow from the assumption
$\Hom(\mathcal B,\mathcal A') = 0$, $\Ext^1(\mathcal B,\mathcal A') =
0$, $\Ext^2(\mathcal C,\mathcal A') = 0$.

We then apply the functor ${\mathbf R}\Hom(\Ker \beta,\bullet)$ to the
exact sequence
\begin{equation*}
   0 \to \mathcal A' \to \Ker \beta' \to E' \to 0
\end{equation*}
to get a short exact sequence
\begin{equation*}
   0 \to \Hom(\Ker \beta,\mathcal A')
   \to \Hom(\Ker \beta,\Ker \beta') \to \Hom(\Ker \beta, E') \to 
   0
   \begin{NB}
     = \Ext^1(\Ker \beta,\mathcal A')
   \end{NB}
   .
\end{equation*}
Therefore the composition of $\varphi$ with the surjection $\Ker
b\twoheadrightarrow E$ can be lifted to a homomorphism
$\psi\colon \Ker \beta\to \Ker \beta'$, unique up to
$\Hom(\Ker \beta,\allowbreak
\mathcal A')$. Therefore we get the claim.

Finally we need to prove that there exists a unique homomorphism
$b\colon \mathcal B\to \mathcal B'$ completing the commutative diagram
(with exact rows)
\begin{equation*}
  \begin{CD}
    0 @>>> \Ker \beta @>>> \mathcal B @>>> \mathcal C @>>> 0
\\
    @. @V{\psi}VV @VV{b}V @VVV @.
\\
    0 @>>> \Ker \beta' @>>> \mathcal B' @>>> \mathcal C' @>>> 0.
  \end{CD}
\end{equation*}
Applying ${\mathbf R}\Hom(\bullet,\mathcal B')$ to the exact sequence
\eqref{eq:exact1}, we get
\begin{equation*}
   \Hom(\mathcal B,\mathcal B') \cong \Hom(\Ker \beta,\mathcal B')
\end{equation*}
from the assumption $\Hom(\mathcal C,\mathcal B') = 0$,
$\Ext^1(\mathcal C,\mathcal B') = 0$.
Therefore the composition of the injection $\Ker \beta'\to \mathcal B'$
and $\psi$ can be extended uniquely to a homomorphism
$b\in \Hom(\mathcal B,\mathcal B')$.
\end{proof}

If both monads $M$, $M'$ are given by complexes of forms
\eqref{eq:cpx} corresponding to data $X$, $X'$, then the vanishing of
five cohomology groups in the assumption of \lemref{lem:monad} can be
easily checked. 
(We only need to compute cohomology groups of line bundles.)
And it is also easy to check
\(
   \Ext^1(\mathcal C,\mathcal A') \cong
   \Hom(V_1, V'_0).
\)
Moreover if we represent a homomorphism $M\to M'$ with respect to the
decomposition in \eqref{eq:cpx}, we get
\begin{Lemma}\label{lem:monadhom}
  A homomorphism $(a,b,c)$ of monads from $M$ to $M'$ is represented
  as
  \begin{gather*}
     a =
     \begin{pmatrix}
       \xi_0 + a_{01} d & a_{01} s
       \\
       0 & \xi_1
     \end{pmatrix},
\quad
     b =
     \begin{pmatrix}
       \id_{\C^2} \otimes (\xi_0 + a_{01} d) 
         & \id_{\C^2} \otimes a_{01} & 0
       \\
       0 & \id_{\C^2} \otimes \xi_1 & 0 
       \\
       0 & 0 & b_{WW}
     \end{pmatrix},
\\
     c =
     \begin{pmatrix}
       \xi_0  & a_{01} s
       \\
       0 & \xi_1 + d' a_{01}
     \end{pmatrix}
  \end{gather*}
satisfying the relations \eqref{eq:commute}, where
\begin{equation*}
  \xi_0 \colon V_0 \to V'_0, \quad
  \xi_1 \colon V_1 \to V'_1, \quad a_{01}\colon V_1\to V'_0, \quad
  b_{WW} \colon W \to W'.
\end{equation*}
\end{Lemma}

The above lemma gives an exact sequence
\begin{equation*}
   0 \to \Hom(V_1, V'_0) \to \Hom_{\mathrm{mon}}(M,M')
   \to \Hom(X,X') \to 0.
\end{equation*}
Moreover any element in $\Hom(V_1,V'_0)$ maps $\Ker\beta$ to
$\Ima\alpha'$. Therefore $\Hom(V_1,V'_0)$ is contained in the kernel
of $H\colon \Hom_{\mathrm{mon}}(M,M')\to \Hom(E,E')$, which we identified with
$\Ext^1(\mathcal C,\mathcal A) \cong \Hom(V_1,V'_0)$ in
\lemref{lem:monad}. Therefore $\Hom(V_1,V'_0)$ is exactly the same as
the kernel of $H$. We conclude
\begin{equation*}
   \Hom(X,X') \cong \Hom(E,E').
\end{equation*}
This completes the proof of \propref{prop:hom}.

\begin{proof}[Proof of \lemref{lem:monadhom}]
We write $a$, $b$, $c$ in matrix forms:
\begin{equation*}
   a =
   \begin{bmatrix}
     a_{00} & a_{01} s \\
     0 & a_{11} 
   \end{bmatrix}, \qquad
   b = 
   \begin{bmatrix}
     b_{00} & b_{01} & b_{0W} \\
     b_{10} & b_{11} & b_{1W} \\
     b_{W0} & b_{W1} & b_{WW}
   \end{bmatrix}
   , \qquad
   c =
   \begin{bmatrix}
     c_{00} & c_{01} s \\
     0 & c_{11} 
   \end{bmatrix},
\end{equation*}
where $b_{ij}\in \Hom(\C^2\otimes V_j,\C^2\otimes V_i)$, etc. Note $a$
and $c$ are upper triangular as $H^0(\shfO(-C)) = 0$. We also used
$H^0(\shfO(C)) \cong \C s$.

We restrict the monads $M$, $M'$ to $\linf = \{ z_0 = 0\}$. Then the
commutativity of the diagram implies
\begin{gather*}
    b_{00} = \operatorname{id}_{\C^2}\otimes a_{00}, \quad
    b_{10} = 0, \quad b_{W0} = 0, \quad
    b_{01} = \operatorname{id}_{\C^2}\otimes a_{01}, 
\\
    b_{11} = \operatorname{id}_{\C^2}\otimes a_{11}, \quad
    b_{W1} = 0,
\\
    b_{00} = (c_{00} + c_{01} d)\otimes\operatorname{id}_{\C^2}, \quad
    d' b_{00} = c_{11} \otimes\operatorname{id}_{\C^2}, \quad
    b_{01} = c_{01} \otimes\operatorname{id}_{\C^2},
\\
    d' b_{01} + b_{11} = c_{11}\otimes\operatorname{id}_{\C^2}, \quad
    b_{0W} = 0, \quad b_{1W} = 0.
\end{gather*}
So linear maps are determined from
$c_{00}$, $a_{11}$, $a_{01}$, $b_{WW}$. We write the first two by
$\xi_0$, $\xi_1$. 

We then consider the commutativity of the diagram on arbitrary point
in $\bp$ to find that $\xi_0$, $\xi_1$ satisfy the relations
\eqref{eq:commute}.
\end{proof}

\subsection{Vanishing of homomorphisms}

Recall that we take $\zeta$ with \eqref{eq:overallass}. Then the
$\zeta$-stability implies (S2) by \lemref{lem:S2}.

\begin{Proposition}
  Suppose that $X = (B_1,B_2,d,i,j)$ is $\zeta$-stable and let $E =
  \Ker\beta/\Ima \alpha$ as above. Then
  \begin{enumerate}
  \item $\Hom(E(-mC),\shfO_C(-1)) = 0$,
  \item $\Hom(\shfO_C,E(-mC)) = 0$.
  \end{enumerate}
\end{Proposition}

\begin{NB}
\begin{Proposition}
  Suppose $(B_1,B_2,d,i,j)$ satisfies the condition \textup{(S2)} in
  \thmref{thm:King} and hence gives a coherent sheaf $E$ on
  $\bp$. Then $(B_1,B_2,d,i,j)$ is $\zeta$-semistable if and only if
\begin{enumerate}
\item If $m\zeta_0 + (m+1)\zeta_1 > 0$, then
$\Hom(\shfO_C(-m-1), E) = 0$.
\item If $m\zeta_0 + (m+1)\zeta_1 < 0$, then
$\Hom(E,\shfO_C(-m-1)) = 0$.
\end{enumerate}
\end{Proposition}
\end{NB}

\begin{proof}
We have $m \zeta_0 + (m+1)\zeta_1 < 0$ and $(m-1)\zeta_0 + m\zeta_1 %
\begin{NB}
  = m\zeta_0 + (m+1)\zeta_1 - (\zeta_0 + \zeta_1)
\end{NB}
> 0$ by the assumption \eqref{eq:overallass}. Hence if $m > 0$, we
have $\Hom(E(-mC),\shfO_C(-1)) = \Hom(E,\shfO_C(-m-1)) = \Hom(X,C_m) =
0$ and $\Hom(\shfO_C,E(-mC)) = \Hom(\shfO_C(-m),E) = \Hom(C_{m-1},X) =
0$ by \propref{prop:ss}. Therefore we only need to show (2) for $m=0$.
(In fact, the following argument only uses that $E$ arises as
$\Ker\beta/\Ima\alpha$. But $m > 0$ case is already proved as
$\Hom(\shfO_C(-m),E) = 0$ implies $\Hom(\shfO_C,E) = 0$.)

Let us take a representation of the quiver in \thmref{thm:King} with
$V_0 = \C$, $V_1 = W = 0$, $B_1 = B_2 = d = 0$. Let us denote it by
$C'$. If we put it in the complex \eqref{eq:cpx}, we get
\begin{equation*}
    \shfO(C-\linf) \xrightarrow{\alpha = \left[
        \begin{smallmatrix}
          z \\ w
        \end{smallmatrix}
        \right]}
      \shfO^{\oplus 2} 
      \xrightarrow{\beta = \left[
          \begin{smallmatrix}
             z_2 & -z_1
          \end{smallmatrix}\right]}
      \shfO(\linf).
\end{equation*}
This is exact except that we have $\operatorname{Cok}\beta = \shfO_C$.
Therefore the complex is quasi-isomorphic to $\shfO_C[-1]$, where
$[-1]$ denote the shift functor of degree $(-1)$.
Let
\( 
  M\colon \mathcal A \xrightarrow{\alpha} \mathcal B
   \xrightarrow{\beta} \mathcal C
\)
be this complex with shifted of degree $1$ (so that it is quasi-isomorphic to 
$\shfO_C$), and 
\(
  M'\colon \mathcal A' \xrightarrow{\alpha'} \mathcal B'
   \xrightarrow{\beta'} \mathcal C'
\)
be the complex \eqref{eq:cpx} corresponding to the data $X$.
This $M$ is not a monad in our sense, but we can define a homomorphism
as a commuting diagram
\begin{equation*}
  \begin{CD}
   \mathcal A @>{\alpha}>> \mathcal B
   @>{\beta}>> \mathcal C @.
\\
   @. @VVbV @VVcV @.
\\
   @. \mathcal A' @>>{\alpha'}> \mathcal B'
   @>>{\beta'}> \mathcal C'.
  \end{CD}
\end{equation*}
Let $\Hom_{\mathrm{mon}}(M,M')$ denote the space of homomorphisms.

\begin{Claim}
\(
   \Hom_{\mathrm{mon}}(M,M') \cong \Hom(\shfO_C,E).
\)  
\end{Claim}

In fact, it is enough to show that we can complete uniquely the following
commutative diagram with exact rows, for any given homomorphism
$\varphi\colon \shfO_C\to E$:
\begin{equation*}
  \begin{CD}
    0 @>>> \Ima\beta @>>> \mathcal C @>>> \shfO_C @>>> 0
\\
    @. @VVV @VVV @VV{\varphi}V @.
\\
    0 @>>> \mathcal A' @>>> \Ker \beta' @>>> E @>>> 0.
  \end{CD}
\end{equation*}
We consider the long exact sequence for the second row tensored with
$\mathcal C^*$. Then we have
\begin{equation*}
  \Hom(\mathcal C,\mathcal A') = 0
\end{equation*}
(this can be checked directly, or follows from $\Hom(\mathcal
C,\mathcal B') = 0$) and
\begin{equation*}
   \Ext^1(\mathcal C,\mathcal A') = 0
\end{equation*}
(this follows from our choice $V_1 = 0$ for the data corresponding to
$\shfO_C$). Thus we have
\begin{equation*}
  \Hom(\mathcal C,\Ker\beta') \cong 
  \Hom(\mathcal C,E).
\end{equation*}
Now the claim follows.

As $\Hom(\mathcal C,\mathcal B') = 0 = \Hom(\mathcal B,\mathcal
A')$, we have $\Hom_{\mathrm{mon}}(M,M') = 0$.
\end{proof}

Now we have proved that the sheaf $E(-mC)$ is perverse coherent, so
completed the half of the proof of \thmref{thm:main}.

\begin{NB}
In the following subsection we check the vanishing of higher direct images
$R^1p_*(E) = 0$.
But this is not necessary now.
 
\subsection{Vanishing of higher direct images}\label{subsec:higherdirect}

\begin{Lemma}
We have $R^1p_*(E(-mC)) = 0$.
\end{Lemma}

\begin{proof}
After tensoring with $\shfO(-mC)$, we apply ${\mathbf R}p_*$ to the
complex \eqref{eq:cpx}.
From the exact sequence
\begin{equation*}
   0 \to \shfO_{\bp}(-mC) \to \shfO_{\bp}((1-m)C) \to \shfO_C(m-1)\to 0
\end{equation*}
and $R^1p_*(\shfO_C(m-1)) = 0$ for $m\ge 0$, we have
$R^1p_*(\shfO_{\bp}(-mC)) \cong R^1p_*(\shfO_{\bp}(1-m)C)$. By the induction
and $R^1p_*\shfO_{\bp} = 0$, we have $R^1p_*(\shfO_{\bp}(-mC)) = 0$ for $m\ge
-1$. Moreover we have 
\begin{gather*}
   p_*(\shfO_{\bp}(C)) \cong p_*(\shfO_{\bp}) \cong \shfO_{\proj^2},
\\
\intertext{and}
\begin{CD}
0 @>>> p_*(\shfO_{\bp}(-mC)) @>>>  p_*(\shfO_{\bp}((1-m)C)) @>>>
   p_*(\shfO_C(m-1)) @>>> 0
\\
   @. @. @. @|
\\
  @. @. @. 
  \makebox[0cm][c]{$\C_0 \otimes H^0(\proj^1,\shfO_{\proj^1}(m-1))$}
\end{CD}
\end{gather*}
from the above exact sequences with $m=0$ and $m\ge 1$ respectively.
From the first isomorphism and the second exact sequence we have
$p_*(\shfO_{\bp}(-mC)) \cong \mathfrak m_0^m$ for $m\ge 1$, where
$\mathfrak m_0$ is the maximal ideal at the origin.
Therefore ${\mathbf R}p_*(E(-mC))$ is equal to
\begin{equation*}
\begin{CD}
\begin{matrix}
   V_0 \otimes \mathfrak m_0^{m-1}\otimes \shfO_{\proj^2}(-\linf)
\\
   \oplus
\\ V_1\otimes \mathfrak m_0^m\otimes \shfO_{\proj^2}(-\linf))
\end{matrix}
@>{p_*(\alpha)}>>
\begin{matrix}
   \C^2\otimes V_0\otimes\mathfrak m_0^m 
   \\ \oplus \\ 
   \C^2\otimes V_1\otimes\mathfrak m_0^m 
   \\
   \oplus \\ W\otimes\mathfrak m_0^m 
\end{matrix}
@>{p_*(\beta)}>>
\begin{matrix}
  V_0\otimes \mathfrak m_0^m\otimes \shfO_{\proj^2}(\linf)
\\
  \oplus
\\
  V_1\otimes \mathfrak m_0^{m+1}\otimes\shfO_{\proj^2}(\linf)
\end{matrix},
\end{CD}
\end{equation*}
for $m\ge 1$ where $p_*(\alpha)$, $p_*(\beta)$ are induced
homomorphisms from $\alpha$, $\beta$ in \eqref{eq:cpx}. For $m=0$, we
have the same complex if we understood $\mathfrak m_0^{-1}$ as
$\shfO_{\proj^2}$.%
\begin{NB2}
Thus the left most term is
\(
\begin{smallmatrix}
   V_0 \otimes\shfO_{\proj^2}(-\linf)
\\
   \oplus
\\ V_1\otimes\otimes \shfO_{\proj^2}(-\linf))
\end{smallmatrix}.
\)
\end{NB2}

Now to show $R^1p_*(E(-mC)) = 0$ it is enough to prove $p_*(\beta)$ is
surjective.%
\begin{NB2}
First consider the case $m=0$.
As $d$ is surjective by \lemref{lem:d_surj} and
$(B_1d,B_2d,i,jd)$ satisfies the stability condition for the framed
moduli space of torsion free sheaves on $\proj^2$, an element of the
form
\(
\left[
\begin{smallmatrix}
  v_0 \\ 0 
\end{smallmatrix}\right]
\)
in the third term is an image of an element of the form
\(
\left[
\begin{smallmatrix}
  v'_0 \\ -d v'_0 \\ w'
\end{smallmatrix}\right]
\)
in the second term.
On the other hand, $\C^2\otimes\shfO_{\proj^2}\to \mathfrak
m_0\otimes\shfO_{\proj^2}(\linf)$ is surjective. Therefore
$p_*(\beta)$ is surjective, hence $R^1p_*(E) = 0$.
\end{NB2}

%
By Nakayama's lemma we can replace the target by
\(
  V_0\otimes \left(\nicefrac{\mathfrak m_0^m}{\mathfrak m_0^{m+1}}\right)
  \oplus
  V_1\otimes \left(\nicefrac{\mathfrak m_0^{m+1}}{\mathfrak m_0^{m+2}}\right)
  .
\)
From the explicit expression of $\beta$ in \eqref{eq:cpx}, we are
reduced to check that 
\begin{equation*}
\begin{CD}
\begin{matrix}
   \C^2\otimes V_0\otimes
   \left(\nicefrac{\mathfrak m_0^m}{\mathfrak m_0^{m+1}}\right)
   \\ \oplus \\ 
   \C^2\otimes V_1\otimes
   \left(\nicefrac{\mathfrak m_0^m}{\mathfrak m_0^{m+1}}\right)
   \\
   \oplus \\ W\otimes
   \left(\nicefrac{\mathfrak m_0^m}{\mathfrak m_0^{m+1}}\right)
\end{matrix}
@>{\beta'}>>
\begin{matrix}
  V_0\otimes \left(\nicefrac{\mathfrak m_0^m}{\mathfrak m_0^{m+1}}\right)
\\
  \oplus
\\
  V_1\otimes \left(\nicefrac{\mathfrak m_0^{m+1}}{\mathfrak m_0^{m+2}}\right)
\end{matrix}
\end{CD}
\end{equation*}
is surjective, where
\begin{equation*}
  \beta' =    \begin{bmatrix}
    0   & 0     & B_2 z_0 & -B_1 z_0 & i z_0\\
    d z_2  & - d z_1 & z_2    & -z_1       & 0
   \end{bmatrix}.
\end{equation*}
If we compose $\beta'$ with the projection to the second factor
\(
  V_1\otimes \left(\nicefrac{\mathfrak m_0^{m+1}}
    {\mathfrak m_0^{m+2}}\right),
\)
then it is clearly surjective.
\begin{NB2}
In fact, we can even restrict the domain to
$\C^2\otimes V_1\otimes
   \left(\nicefrac{\mathfrak m_0^m}{\mathfrak m_0^{m+1}}\right)$.
\end{NB2}
Therefore it is enough to show that the subspace
\(
  V_0\otimes \left(\nicefrac{\mathfrak m_0^m}{\mathfrak
    m_0^{m+1}}\right)
  \oplus 0
\)
of the target, consisting elements with vanishing the second component,
is contained in the image.
We restrict the domain to the subspace consisting of elements of the
form
\begin{gather*}
\begin{bmatrix}
   v'_0
\\ 
   -dv_0 + \left[
     \begin{smallmatrix}
       z_1 \\ z_2
     \end{smallmatrix}\right]
   v'_1  
\\ 
   w'
\end{bmatrix}
\in
\begin{matrix}
   \C^2\otimes V_0\otimes
   \left(\nicefrac{\mathfrak m_0^m}{\mathfrak m_0^{m+1}}\right)
   \\ \oplus \\ 
   \C^2\otimes V_1\otimes
   \left(\nicefrac{\mathfrak m_0^m}{\mathfrak m_0^{m+1}}\right),
   \\
   \oplus \\ W\otimes
   \left(\nicefrac{\mathfrak m_0^m}{\mathfrak m_0^{m+1}}\right)
\end{matrix}
\\
\intertext{where}
   v'_0 \in\C^2\otimes V_0\otimes
   \left(\nicefrac{\mathfrak m_0^m}{\mathfrak m_0^{m+1}}\right),
\qquad
   v'_1 \in V_1\otimes
   \left(\nicefrac{\mathfrak m_0^{m-1}}{\mathfrak m_0^{m}}\right), 
\qquad
   w'\in W\otimes
   \left(\nicefrac{\mathfrak m_0^m}{\mathfrak m_0^{m+1}}\right),
\end{gather*}
where $\nicefrac{\mathfrak m_0^{m-1}}{\mathfrak m_0^{m}}$ is
understood as $0$ when $m=0$.
Then the restriction is the same as $\tau$ in \eqref{eq:cpx+},
therefore it is surjective by \lemref{lem:Ext+}(3).
\end{proof}
\end{NB}

\section{Inverse construction}\label{sec:inverse}

We give the inverse construction from a sheaf to ADHM data.

\subsection{Two spectral sequences from the resolution of the diagonal}

Let $p_1$, $p_2\colon\bp\times\bp\to \bp$ be the projections to the
first and second factors respectively. Let $C^\bullet$ be the
resolution of the diagonal in \eqref{eq:diagonal}.
Suppose that $(E,\Phi)$ is a framed sheaf on (the first factor) $\bp$
such that $(E(-mC),\Phi)$ is stable perverse coherent. We consider the
double complex for the hyper direct image
\(
  R^\bullet p_{2*}(p_1^* E(-\linf)\otimes C^\bullet)
\)
for which we can consider two spectral sequences as usual. If we take
the cohomology of $C^\bullet$ first, we get 
\(
  E_2^{p,q} = E(-\linf)
\)
for $(p,q)=(0,0)$ and $=0$ otherwise. Thus the spectral sequence
degenerates at $E_2$-term and converges to $E(-\linf)$. 

On the other hand, if we take the direct image first, we get
\linebreak[4]
\(
  R^q p_{2*} (p_1^* E(-\linf)\otimes C^p).
\)
As each $C^p$ is of a form $C^p = p_1^*\mathcal F_1\otimes
p_2^*\mathcal F_2$, this is equal to
\(
  \mathcal F_2\otimes H^q(\bp,E(-\linf)\otimes \mathcal F_1).
\)

\subsection{Vanishing of cohomology groups}

\begin{Lemma}\label{lem:cohvanish}
Suppose that $(E(-mC),\Phi)$ is a stable perverse coherent framed
sheaf with $m\ge 0$. Then

\textup{(1)} $H^q(\bp,E(-k\linf)) = 0$, 
$H^q(\bp,E(C-k\linf)) = 0$ for $k=1,2$, $q=0,2$.

\textup{(2)} $H^q(\bp,E(-\linf)\otimes\mathcal Q^\vee(C)) = 0$
for $q=0,2$.
\end{Lemma}

\begin{proof}
We take direct images of $E$ and $E(C)$ with respect to
$p\colon \bp\to\proj^2$ and deduce the assertion from the vanishing
theorem for sheaves on $\proj^2$ in \cite[Lemma~2.4]{Lecture}.

As $E(-mC)$ is stable perverse, we have $0 = \Hom(\shfO_C,E(-mC)) =
\Hom(\shfO_C(-m),E)$ by (3) in the introducion. Therefore
$\Hom(\shfO_C, E) = 0$ as we have
$\Hom(\shfO_C(-m),\shfO_C)\otimes\shfO_C(-m)\to \shfO_C$ is
surjective. Similarly we have $\Hom(\shfO_C,E(C)) = \Hom(\shfO_C(1),E)
= 0$.%
\begin{NB}
If $\Hom(\shfO_C(-m),E) = 0$, then $\Hom(\shfO_C(n),E) = 0$ for all
$n\ge -m$.
\end{NB}
Note that the condition (3) is equivalent to the torsion freeness of
$p_*E$ at $[1:0:0]$.
In fact, if $\C_0$ is the skyscraper sheaf at $[1:0:0]$ of $\proj^2$,
  then we have
$\Hom(\shfO_C,E) = \Hom(p^*\C_0, E) = \Hom(\C_0,p_*E)$.
Therefore both $p_*E$ and $p_*(E(C))$ are torsion-free in our
situation as they are also torsion free on $\bp\setminus C
\xrightarrow[\cong]{p} \proj^2\setminus\{[1:0:0]\}$ by the condition
(4)'.
Therefore we have 
\(
  H^0(\bp,E(-k\linf)) \cong H^0(\proj^2,(p_*E)(-k\linf)) = 0
\)
and
\(
  H^0(\bp,E(C-k\linf)) \cong H^0(\proj^2,(p_* E(C))(-k\linf)) = 0
\)
for $k=1,2$ by \cite[Lemma~2.4]{Lecture}.

Let us turn to $H^2$. As $R^1 p_*(E)$ is supported at the origin, it
does not contribute to $H^2$. Therefore we have $H^2(\bp,E(-k\linf))
\cong H^2(\proj^2,(p_*E)(-k\linf))$. Hence we have $H^2(\bp,E(-k\linf))
= 0$ for $k=1,2$ by \cite[Lemma~2.4]{Lecture}. The same argument works
for $E(C)$.

For (2) we recall that $\Qcal$ is a pull-back of a vector bundle
$\bar\Qcal$ over $\proj^2$. Its definition is exactly the same as the
vector bundle $\shfO_{\proj^2}^{\oplus 3}/\shfO_{\proj^2}(-\linf)$
appeared in \cite[\S2.1]{Lecture}. Therefore the assertion follows
from \cite[Lemma~2.4]{Lecture}.
\end{proof}

From this vanishing theorem, the spectral sequence must degenerate at the $E_2$-term. Furthermore, since the spectral sequence must converges to
$E(-\linf)$, we find that (1) $E$ is isomorphic to
$\Ker\beta/\Ima\alpha$ where
\begin{equation*}
  \begin{CD}
\begin{matrix}
   V_0 \otimes \shfO(C-\linf)
\\
   \oplus
\\ V_1\otimes\shfO(-\linf)
\end{matrix}
@>{\alpha}>>
\begin{matrix}
   \C^2\otimes V_0\otimes \shfO
   \\ \oplus \\ 
   \widetilde W\otimes \shfO
\end{matrix}
@>{\beta}>>
\begin{matrix}
  V_0\otimes\shfO(\linf) \\ \oplus \\
  V_1\otimes\shfO(-C+\linf)
\end{matrix},
\end{CD}
\end{equation*}
where
\begin{equation*}
  \begin{gathered}
   V_0 = H^1(E(-\linf)), \quad
   V_1 = H^1(E(C-2\linf)), \\
   \widetilde W = H^1(E(-\linf)\otimes \Qcal^\vee(C)),
  \end{gathered}
\end{equation*}
and (2) $\alpha$ is injective and $\beta$ is surjective.  When we have
identified the term of the complex, we have used that the
multiplication by $z_0$ induces isomorphisms $H^1(E(C-2\linf))\cong
H^1(E(C-\linf))$, which is proved as in \cite[Lemma~2.4]{Lecture}
during the course of the above proof.

\subsection{A normal form}
We next show that the above complex is of the form in \eqref{eq:cpx}.
The homomorphisms $\alpha$, $\beta$ are induced from $d^{-1}$, $d^0$
in \eqref{eq:diagonal}. Hence we have
\begin{equation*}
  \begin{split}
    &
  \alpha = 
     \begin{bmatrix}
    \id\otimes z & B_1\otimes z_0\\
    \id\otimes w & B_2\otimes z_0\\
    0 & \alpha_0 z_0 + \alpha_1 z_1 + \alpha_2 z_2
   \end{bmatrix},
\\
   & \beta =
   \begin{bmatrix}
     \id\otimes z_2 & - \id\otimes z_1 & 
     \beta_0 z_0 + \beta_1 z_1 +\beta_2 z_2
     \\
     d \otimes w & - d \otimes z & 
     \gamma_0 w + \gamma_1 z
   \end{bmatrix}
  \end{split}
\end{equation*}
for $\alpha_0$, $\alpha_1$, $\alpha_2\in \Hom(V_1,\widetilde W)$,
$\beta_0$, $\beta_1$, $\beta_2\in \Hom(\widetilde W,V_0)$,
$\gamma_0$, $\gamma_1\in \Hom(\widetilde W,V_1)$. Also $B_1$, $B_2$,
$d$ are induced homomorphisms of $z'$, $w'$, $s'$ respectively.

From the form of $d^{0}$ and $\chi$, we have
$\beta_1 = \beta_2 = 0$. Now we restrict the complex to $\linf$ to get
\begin{equation*}
  V_1\otimes \shfO_{\linf}(-\linf)
    \xrightarrow{\alpha_1 z_1 + \alpha_2 z_2}
  \widetilde W\otimes \shfO_\linf
    \xrightarrow{\gamma_0 w + \gamma_1 z}
  V_1\otimes \shfO_\linf(\linf).
\end{equation*}
Then the argument in \cite[\S2]{Lecture} shows that the trivialization
$E|_\linf\to \shfO_\linf^{\oplus r}$ induces a decomposition
$\widetilde W\to V_1\oplus V_1 \oplus W$ such that
\begin{equation*}
   \alpha_1 =
   \begin{bmatrix}
     \id_{V_1} \\ 0 \\ 0
   \end{bmatrix},
   \ 
   \alpha_2 =
   \begin{bmatrix}
     0 \\ \id_{V_1} \\ 0
   \end{bmatrix},
   \ 
   \gamma_0 =
   \begin{bmatrix}
     \id_{V_1} & 0 & 0
   \end{bmatrix},
      \ 
   \gamma_1 =
   \begin{bmatrix}
     0 & - \id_{V_1} & 0
   \end{bmatrix}
\end{equation*}
From the equation $\beta\alpha = 0$, we then get
\begin{equation*}
   \alpha_0 =
   \begin{bmatrix}
      - d B_1 \\ - d B_2 \\ j
   \end{bmatrix},
   \quad
   \beta_0 =
   \begin{bmatrix}
      B_2 & -B_1 & i
   \end{bmatrix}
\end{equation*}
for $i\colon W\to V_0$, $j\colon V_1\to W$.

Thus the complex is of the form in \eqref{eq:cpx}.

\subsection{$\zeta$-stability}
By \lemref{lem:King} the data $X = (B_1,B_2,d,i,j)$ satisfies the
condition (S2) in \thmref{thm:King}.
As $E(-mC)$ is perverse coherent, we have
\begin{equation*}
\begin{split}
   & \Hom(X,C_n) \cong \Hom(E(-mC),\shfO_C(m-n-1)) = 0,
\\  
   & \Hom(C_{n'},X) \cong \Hom(\shfO_C(m-n'-1), E(-mC)) = 0
\end{split}
\end{equation*}
for $n' < m \le n$. By \propref{prop:ss2} $X$ is
$\zeta$-stable for $\zeta$ satisfying \eqref{eq:overallass}. 

\subsection{The last part of the proof}

Let us finally comment that maps constructed in the
previous section and this section are mutually inverse.

The statement `sheaf $\to$ ADHM $\to$ sheaf is an identity' is proved
already in our construction of the data $X = (B_1,B_2,d,i,j)$ via the
resolution of the diagonal.

The statement `ADHM $\to$ sheaf $\to$ ADHM is an identity' means the
uniqueness of $X$ for a given framed sheaf $E$. But this is also
already proved in \propref{prop:hom}.

\begin{Remark}
  As a corollary of our proof, we have an equivalence between the
  derived category $D^b(\operatorname{Coh}_{C}(\bp))$ of complexes of
  coherent sheaves whose homologies are supported on $C$
  set-theoretically, and the
  derived category of finite dimensional nilpotent modules of the quiver 
\(
\xymatrix{
  \bullet \ar@<-1ex>[r]_{d} & \ar@<-1ex>[l]_{B_1,B_2} \ar[l]
  \bullet
}
\)
with relation $B_1 d B_2 = B_2 d B_1$. In fact, the proof is much
simpler, as we do not need to study the stability condition.
\end{Remark}

\section{Distinguished chambers}

\subsection{Torsion free sheaves on blowup}\label{subsec:dist}

There are two distinguished chambers in the space of parameters for
the stability conditions. The first one is the case when the stability
condition is given by (S1)' and (S2). This is the chamber adjacent to
the boundary $\zeta_0+\zeta_1 = 0$ of the region \eqref{eq:ass}. Let
us take a parameter ${}^\infty\zeta$ from the chamber as in
Figure~\ref{fig:zeta1}.
The corresponding moduli space is the framed moduli space of torsion
free sheaves on $\bp$. This follows from
\begin{Proposition}
Fix $r$, $k$, $n$. If $m$ is sufficiently large, $E(-mC)$ is
stable perverse coherent if and only if $E$ is torsion free.
\end{Proposition}

This, in turn, follows from the following lemma once we observe that
$X$ is ${}^\infty\zeta$-stable if and only if it satisfies the
condition (S2) and (S1)',
\begin{NB}
  take $\zeta = (-1,m/m+1)$ and consider the conditions in
  \defref{def:stable} in the limit $m\to\infty$.
\end{NB}%
and the sheaf $E$ corresponding to $X$ is torsion free if and only
$\alpha$ is injective except finitely many points.
We prove \thmref{thm:King} simultaneously. For this we use the
observation that the sheaf $E$ corresponding to $X$ is locally free if
and only $\alpha$ is injective at every point in $\bp$.

\begin{Lemma}
Let $X = (B_1,B_2,d,i,j)\in\mu^{-1}(0)$. Then

\textup{(1)(\cite[4.1.3]{King})} $X$ satisfies the condition \textup{(S1)} in
\thmref{thm:King} if and only if $\alpha$ in the complex \eqref{eq:cpx}
is injective at every point in $\bp$.

\textup{(2)} $X$ satisfies the condition \textup{(S1)'} in
\secref{sec:ADHM} if and only if $\alpha$ is injective possibly except
finitely many points.
\end{Lemma}

\begin{proof}
The proof is taken from that in \cite[Prop.~4.1]{Na:ADHM}. In fact, it
becomes much simpler.

Let us first show the `only if' part: Assume $\alpha$ is not injective
at a point $([z_0:z_1:z_2],[z:w])$. As $\alpha$ is injective along
$z_0=0$, we may assume $z_0 = 1$. Then there exists $v_0\in V_0$,
$v_1\in V_1$ such that
\begin{equation*}
  \begin{gathered}
  (z_1 - d B_1) v_1 = 0, \quad
  (z_2 - d B_2) v_1 = 0, \\
  B_1 v_1 + z v_0 = 0, \quad
  B_2 v_1 + w v_0 = 0, \quad
  j  v_1 = 0
  \end{gathered}
\end{equation*}
and $(v_0,v_1)\neq (0,0)$. From the equation, we have $v_1\neq 0$, as
$v_1 = 0$ implies $v_0 = 0$.
Then $(S_0,S_1) = (\C v_0,\C v_1)$ violates the condition (S1). This
completes the proof of `only if' part for (1).

By the condition (S1)' we must have $\dim S_0 = 1$, i.e., $v_0\neq
0$. We consider $X'$ induced (from $X$) on the quotient $V'_0 :=
V_0/S_0$, $V'_1 := V_1/S_1$. 

We claim $X'$ satisfies the condition (S1)'. Let $S'_0\subset V'_0$,
$S'_1\subset V'_1$ be subspaces as in the condition (S1)'. Then their
inverse images $\widetilde S'_0$, $\widetilde S'_1$ in $V_0$, $V_1$
must satisfy $\dim \widetilde S'_0 \ge \dim \widetilde S'_1$ or
$\widetilde S'_0 = \widetilde S'_1 = 0$ as $X$ satisfies the condition
(S1)'. As $\dim S_0 = \dim S_1 = 1$, we have $\dim S'_0 \ge \dim
S'_1$. This completes the proof of the claim.

By the induction on the dimensions of $V_0$, $V_1$ we may assume that
$\alpha$ defined for $X'$ is injective possibly except for finitely
many points.

Let us consider the complex \eqref{eq:cpx} for $X$, $X'$ and
$(S_0,S_1)$ with $(B_1,B_2,d) = (z,w,s)$ as in
\thmref{thm:blowupplane}. We have a commutative diagram:
\begin{equation*}
  \begin{CD}
     C^{-1}(S) @>\alpha>> C^0(S) @>\beta>> C^1(S)
\\
     @VVV @VVV @VVV
\\
     C^{-1}(X) @>\alpha>> C^0(X) @>\beta>> C^1(X)
\\
     @VVV @VVV @VVV
\\
     C^{-1}(X') @>\alpha>> C^0(X') @>\beta>> C^1(X')
  \end{CD}
\end{equation*}
From \lemref{lem:diag}(2) (or a direct argument) $\alpha$ for $S$ is
injective except at the point $([z_0:z_1:z_2],[z:w])$. Therefore
the injectivity of $\alpha$ for $X$ fails only at
this point or a point where $\alpha$ for $X'$ is not
injective. So we only have finitely many points.

The proof of the `if' part can be given exactly as in
\lemref{lem:King}, so is omitted.%
\begin{NB}
We prove the `if' part: Assume $X$ violates the condition (S1)'. We
have subspaces $S_0$, $S_1$ satisfying
$B_\alpha(S_1)\subset S_0$ ($\alpha=1,2$), $d(S_0) \subset S_1$, $\Ker
j\supset S_1$, $\dim S_0 < \dim S_1$. We take $S_0$, $S_1$ among such
pairs of subspaces so that
\begin{equation*}
   \theta(S_0,S_1) := \frac{-\dim S_0 + \dim S_1}{\dim S_0 + \dim S_1}
\end{equation*}
takes the maximum. In particular, $\theta(S_0,S_1) > 0$. Moreover, we
may assume $(B_1,B_2,d)|_{(S_0,S_1)}$ is $\theta$-stable. If we define a
new parameter
\begin{equation*}
   (\zeta'_0,\zeta'_1) = (-1,1) - \theta(S_0,S_1)(1,1),
\end{equation*}
we have $\zeta'_0\dim S_0 + \zeta'_1\dim S_1 = 0$ and the
$\theta$-stability is equivalent to the
$(\zeta'_0,\zeta'_1)$-stability. We moreover have $\zeta'_0+\zeta'_1 <
0$, $\zeta'_0 < 0$. Therefore $(B_1,B_2,d)|_{(S_0,S_1)}$ must be the
one classified in \thmref{thm:W=0}. Then the injectivity of $\alpha$
is violated along the exceptional divisor $C$.

If $X$ violates the condition (S1), we repeat the same argument except
we only have $\theta(S_0,S_1) \ge 0$ and hence $\zeta'_0+\zeta'_1 \le
0$. Then we have other possibility that $(B_1,B_2,d)|_{(S_0,S_1)}$
corresponds to a point in $\bp$. So the injectivity of $\alpha$ is
violated at the point.
\end{NB}
\end{proof}

\thmref{thm:King} with \thmref{thm:blowupplane} implies
$\bM_\zeta^{\mathrm{ss}}(r,k,n)$ is the Uhlenbeck compactification of
the framed moduli $\bM_\zeta^{\mathrm{s}}(r,k,n) = \bMreg(r,k,n)$ of
locally free sheaves on $\bp$ if $\zeta_0+\zeta_1 = 0$, $\zeta_0 < 0$.

\subsection{Another distinguished chamber 
-- torsion free sheaves on blow-down}\label{subsec:another2}

We describe the other distinguished chamber. It is the chamber
$\{ \zeta_0 < 0, \zeta_1 < 0\}$. Let us denote a parameter from this
chamber by $\tzeta$.
In this case the $\tzeta$-stability condition just means that a
pair of subspaces $T_0 \subset V_0$, $T_1 \subset V_1$ such that
$B_\alpha(T_1)\subset T_0$, $d(T_0)\subset T_1$, $\Ima i\subset T_0$
must be $T_0 = V_0$, $T_1 = V_1$. Then both $(dB_1,dB_2,di,j)$ and
$(B_1d,B_2d,i,jd)$ satisfy the stability condition for the framed
moduli space of torsion free sheaves on $\proj^2$ described in the
beginning of \secref{sec:ADHM}. Therefore we have the morphism
$\bM_{\tzeta}(r,k,n)\to M(r,\dim V_0)$ and $\to M(r,\dim V_1)$. This
gives a sharp contrast to moduli spaces from other chambers. In
general, we only have morphisms into the Uhlenbeck (partial)
compactification $M_0(r,\dim V_0)$ or $M_0(r,\dim V_1)$.

\begin{Lemma}\label{lem:d_surj}
  If $[(B_1,B_2,d,i,j)]\in\bM_{\tzeta}(r,k,n)$, then $d$ is
  surjective. In particular, $\bM_{\tzeta}(r,k,n) = \emptyset$ if $k >
  0$.
\end{Lemma}

\begin{proof}
  We take $T_0 = V_0$, $T_1 = \Ima d$. Then all conditions are
  satisfied, so we must have $T_1 = V_1$, i.e., $d$ must be
  surjective.
\end{proof}

We now consider the case $k=0$. Then the situation is surprisingly very
simple:

\begin{Proposition}\label{prop:c_1=0}
  The morphism $\bM_{\tzeta}(r,0,n)\to M(r,n)$ is an isomorphism.
\end{Proposition}

\begin{proof}
  From the assumption $k=0$ and \lemref{lem:d_surj}, $d$ must be an
  isomorphism. Therefore $(B_1,B_2,d,i,j)$ can be recovered from
  \break
  $(dB_1,dB_2, di, j)$. Conversely if $(\tilde B_1,\tilde B_2,\tilde
  i,\tilde j)$ is stable for $M(r,n)$, then $(\tilde B_1, \tilde B_2,
  \operatorname{id}_V, \tilde i,\tilde j)$ gives a $\tzeta$-stable
  representation. It gives the inverse of $M_{\tzeta}(r,0,n)\to M(r,n)$.
\end{proof}

The $c_1\neq 0$ case will be studied in the subsequent paper.

\begin{Remark}\label{rem:zeta1}
  If we take a parameter $\zeta$ from the adjacent chamber, i.e.,
  $\zeta_0 < 0$, $0 < \zeta_1 \ll 1$, the $\zeta$-stability implies
  that a pair of subspaces $T_0 \subset V_0$, $T_1 \subset V_1$ such
  that $B_\alpha(T_1)\subset T_0$, $d(T_0)\subset T_1$, $\Ima i\subset
  T_0$ must be $T_0 = V_0$.
Then $(B_1 d, B_2 d, i, jd)$ satisfy the stability condition for the
framed moduli space of torsion free sheaves on $\proj^2$. Therefore we
have a morphism $\bM_{\zeta}(r,k,n)\to M(r,\dim V_0)$.
\end{Remark}

\end{document}